\documentclass[12pt]{article}
\usepackage{amsmath}
\usepackage{amssymb}
\usepackage{amsmath,amssymb,amsthm,amsfonts}
{\theoremstyle{definition}

}

\newcommand{\il}[2]{\int\limits_{#1}^{#2}}

\newcommand{\ph}{\phantom{a}}
\newcommand{\phh}{\phantom{aaa}}

\newcommand{\pinf}{+\infty}

\newcommand{\no}{\eqno}

\newcommand{\ta}{\tau}
\newcommand{\noin}{\noindent}
\newcommand{\vsk} {\vskip 10pt}

\newcommand{\fr}{\frac}
\newcommand{\sist}[2]{\left\{
\begin{array}{l}
{#1}\\
\ph\\
{#2}
\end{array}
\right.}

\begin{document}

MSC 34D20

\vskip 20pt

\centerline{\bf Comparison criteria for first order polynomial differential equations}

\vskip 10 pt

\centerline{\bf G. A. Grigorian}

\vskip 10 pt

\centerline{0019 Armenia c. Yerevan, str. M. Bagramian 24/5}
\centerline{Institute of Mathematics of NAS of Armenia}
\centerline{E - mail: mathphys2@instmath.sci.am, \ph phone: 098 62 03 05, \ph 010 35 48 61}

\vskip 20 pt

\noindent
Abstract. In this paper we use the comparison method for investigation of  first order polynomial differential equations. We prove two  comparison criteria for these equations. The proved  criteria we use  to obtain some global solvability criteria for first order polynomial differential equations. On the basis od these criteria  we prove some criteria for existence of a closed solution (of closed solutions)   for first order polynomial differential equations. The results obtained we compare with some known results.

\vskip 20 pt


\vskip 20 pt

{\bf  1. Introduction}.  Let $a_k(t), \ph k=\overline{1,n}$  be real-valued continuous functions on $[t_0,\tau_0) \ph (t_0 < \tau_0 \le \pinf)$. Consider the first order polynomial differential equation
$$
y' +  \sum\limits_{n=0}^{n}a_k(t) y^k = 0,  \phh t_0\le t \le \tau_0. \no (1.1)
$$
According to the general theory of normal systems of differential equations for every $t_1  \ge t_0, \ph \gamma \in \mathbb{R}$ and for any solution $y(t)$ of Eq. (1.1) with $y(t_1) = \gamma$ there exists $t_2 > t_2$ such that $y(t)$ is continuable on $[t_1,t_2)$. From the point of view of qualitative theory of differential equations an important interest represents the case $t_2 =+\infty$. One of effective ways to study the conditions, under which the case $t_2 = +\infty$ holds, is the comparison method. This method has been used in [1,2] to obtain some comparison criteria for Eq. (1.1) in the case  $n=2$ (the case of Riccati equations), which were used for qualitative study of different types of equations (see e. g. [3--13]).
In the general case Eq. (1.1) attracts the attention of mathematicians in the connection with a relation of the problem of existence of closed solutions of Eq. (1.1) with the problem of determination of the upper bound for the number of limit cycles in two-dimensional polynomial vector fields of degree n.
(see [] the 16th problem of Hilbert [recall that a solution $y(t)$ of Eq. (1.1), existing on any interval $[t_0,T]$, is called closed on that interval, if $y(t_0) = y(T)$]) and many works are devoted to it (see [15--17] and cited works therein). Significant results in this direction have been obtained in [18]. Among them we point out the following result.

{\bf Theorem 1.1. ([18, p.3, Theorem 1])}. {\it Let us assume that $a_0(t)\equiv 0, \linebreak \il{t_0}{T} a_1(t) d t >~ 0$. Let us assume that there exists some $j=2,\ldots, n$such that $(-1)^k a_k(t) \ge 0$ and $\sum\limits_{k=j}^n(-1)^k a_k(t) > 0$ for all $k=j,\ldots, n$ and $t \in [t_0,T]$. Then there exists a negative isolated closed solution of Eq. (1.1) on $[t_0,T]$.
}

\phantom{aaaaaaaaaaaaaaaaaaaaaaaaaaaaaaaaaaaaaaaaaaaaaaaaaaaaaaaa} $\blacksquare$

This and other theorems of work [18] were obtained by the use of a perturbation method and the contracting mapping principle. In this paper we use the comparison method for investigation of Eq. (1.1) for the case $n \ge 3$. In section 3 we prove two comparison criteria for Eq. (1.1). These criteria we use in section 4 to obtain some global solvability criteria for Eq. (1.1). On the basis od these criteria in section 5 we prove some criteria for existence of a closed solution (of closed solutions) of Eq. (1.1). The results obtained we compare with results of work [18].

\vskip 20 pt

{\bf  2. Auxiliary propositions}. Denote

$$
D(t,u,v) \equiv \sum\limits_{k=1}^n a_k(t) S_k(u,v),
$$
where $S_k(u,v)\equiv \sum\limits_{j=0}^{k-1}u^j v^{k-j-1}, \ph u,v \in \mathbb{R}, \ph k=\overline{1,n}, \ph t \ge t_0.$ Let $b_k(t),\ph k=\overline{0,n}$ be real-valued continuous functions on $[t_0, +\infty)$. Consider the equation
$$
y' + \sum\limits_{k=0}^n b_k(t) y^k = 0, \phh t \ge t_0. \eqno (2.1)
$$
Let $y_0(t)$ and $y_1(t)$ be solutions of the equations (1.1) and (2.1) respectively on $[t_1,t_2)\subset [t_0,+\infty)$. Then
$$
[y_0(t) - y_1(t)]' + \sum\limits_{k=0}^n a_k(t)[y_0^k(t) - y_1^k(t)] + \sum\limits_{k=0}^n[a_k(t) - b_k(t)] y_1^k(t) = 0, \ph t \in [t_1,t_2).
$$
It follows from here and from the obvious equalities $y_0^k(t) - y_1^k(t) = [y_0(t) - y_1(t)] S_k(y_0(t),y_1(t)), \linebreak k= \overline{1,n},$ that
$$
[y_0(t) - y_1(t)]' + D(t,y_0(t),y_1(t))[y_0(t) - y_1(t)] + \sum\limits_{k=0}^n[a_k(t) - b_k(t)]y_1^k(t) = 0, \ph t \in [t_1,t_2).
$$
It is clear from here that $y_0(t) - y_1(t)$ is a solution of the linear equation
$$
x' + D(t,y_0(t),y_1(t))x + \sum\limits_{k=0}^n[a_k(t) - b_k(t)]y_1^k(t) = 0, \ph t \in [t_1,t_2).
$$
Then by the Cauchy formula we have
$$
y_0(t) - y_1(t) = \exp\biggl\{-\il{t_1}{t}D(\tau,y_0(\tau),y_1(\tau)) d \tau\biggr\}\biggl[y_0(t_1) - y_1(t_1) - \phantom{aaaaaaaaaaaaaaaaaaaaaaaaaaaaaaaaaaaaa}
$$
$$
\phantom{aaaaaaa}- \il{t_1}{t}\exp\biggl\{\il{t_1}{\tau}D(s,y_0(s),y_1(s))d s\biggr\}\biggl(\sum\limits_{k=0}^{n}[a_k(\tau) - b_k(\tau)]\biggr) d \tau\biggr], \phh t \in [t_1,t_2). \eqno (2.2)
$$
Consider the differential inequality
$$
\eta' + \sum\limits_{k=0}^n a_k(t) \eta^k \ge 0, \phh t_0 \le t <\ta_0. \no (2.3)
$$

{\bf Definition 2.1.} {\it A continuous on $[t_0,\tau_0) \ph (\tau_0 \le +\infty)$ function $\eta^*(t)$ is called a sub solution of the inequality (2.3) on $[t_0,\tau_0)$ if for every $t_1 \in [t_0,\tau_0)$ there exists a solution $\eta_{t_1}(t)$ of the inequality (2.3) on $[t_0,t_1]$ such that $\eta_{t_1}(t_0) \ge \eta^*(t_0), \ph \eta_{t_1}(t_1) = \eta^*(t_1)$.
}

Consider the differential inequality
$$
\zeta' + \sum\limits_{k=0}^n a_k(t) \zeta^k \le 0, \phh t_0 \le t <\ta_0. \no (2.4)
$$
{\bf Definition 2.2.} {\it A continuous on $[t_0,\tau_0) \ph (\tau_0 \le +\infty)$ function $\zeta^*(t)$ is called a super solution of the inequality (2.4) on $[t_0,\tau_0)$ if for every $t_1 \in [t_0,\tau_0)$ there exists a solution $\zeta_{t_1}(t)$ of the inequality (2.4) on $[t_0,t_1]$ such that $\zeta_{t_1}(t_0) \le \zeta^*(t_0), \ph \zeta_{t_1}(t_1) = \zeta^*(t_1)$.

Obviously any solution $\eta(t) \ph (\zeta(t))$ of the inequality (2.3) ((2.4)) on $[t_0,\tau_0)$ is also a sub (super)  solution of that inequality on $[t_0,\tau_0)$.

\vsk

{\bf Lemma 2.1.} {\it Let $y(t)$ be a solution of Eq. (1.1) on $[t_0,\tau_0)$ and Let $\eta^*(t)$ be a sub solution of the inequality (2.3) on $[t_0,\tau_0)$ such that $y(t_0) \le \eta^*(t_0)$. Then $y(t) \le \eta^*(t), \ph t \in [t_0,\tau_0)$, and if $y(t_0) < \eta^*(t_0)$, then $y(t) < \eta^*(t), \ph t \in [t_0,t_1)$.
}

Proof. It is enough to show that if  $\eta(t)$ is a solution  of the inequality (2.3) on $[t_0,\tau_0)$ with $y(t_0) \le \eta(t_0)$, then
$$
y(t) \le \eta(t), \ph t \in [t_9,\tau_0). \eqno (2.5)
$$
and if $y(t_0) < \eta(t_0)$, then
$$
y(t) < \eta(t), \ph t \in [t_0,\tau_0). \eqno (2.6)
$$
We set $\widetilde{a}_0(t) \equiv -\eta'(t) - \sum\limits_{k=1}^n a_k(t) \eta^k(t), \ph t \in [t_0,\tau_0)$.  By (2.3) we have
$$
\widetilde{a}_0(t) \le a_0(t), \ph t \in [t_0,\tau_0). \eqno (2.7)
$$
Obviously $\eta(t)$ is a solution of the equation
$$
y' + \sum\limits_{k=1}^n a_k(t) y^k + \widetilde{a}_0(t) = 0, \ph t \in [t_0,\tau_0)
$$
on $[t_0,\tau_0)$. Then in virtue of (2.2) we have
$$
y_0(t) - \eta(t) = \exp\biggl\{-\il{t_0}{t}D(\tau,y(\tau),\eta(\tau)) d \tau\biggr\}\biggl[y(t_0) - \eta(t_0) + \phantom{aaaaaaaaaaaaaaaaaaaaaaaaaaaaaaaaaaaaa}
$$
$$
\phantom{aaaaaaaaaaaaaaaaaaa}+ \il{t_1}{t}\exp\biggl\{\il{t_1}{\tau}D(s,y(s),\eta(s))d s\biggr\}\biggl(\widetilde{a}_0(\tau) - a_0(\tau)\biggr) d \tau\biggr], \phh t \in [t_0,\tau_0).
$$
This together with (2.7) implies that if $y(t_0) \le \eta(t_0) \ph (y(t_0) < \eta(t_0))$, then  (2.6) ((2.7)) is valid. The lemma is proved.

By analogy with the proof of Lemma 2.1 one can prove the following lemma

{\bf Lemma 2.2.} {\it Let $y(t)$ be a solution of Eq. (1.1) on $[t_0,\tau_0)$ and let $\zeta^*(t)$ be a super solution of the inequality (2.4) on $[t_0,\tau_0)$ such that $\zeta^*(t_0) \le y(t_0)$. Then $\zeta^*(t) \le y(t), \ph t \in [t_0,\tau_0)$, and if $\zeta(t_0) < y(t_0)$, then $\zeta(t) < y(t), \ph t \in [t_0,\tau_0).$}

\phantom{aaaaaaaaaaaaaaaaaaaaaaaaaaaaaaaaaaaaaaaaaaaaaaaaaaaaaaaaaaaaa}  $\blacksquare$

{\bf Remark 2.1.} {\it It is clear that Lemma 2.1 (Lemma 2.2) remains valid if in the case $\tau_0 < +\infty$ the interval $[t_0,\tau_0)$ is replaced by $[t_0,\tau_0]$ in it.}

For any $T > t_0$ and $j=2,\ldots, n$ we set
$$
M_{T,j}\equiv \max\biggl\{1, \max\limits_{\tau\in[t_0,T]}\biggl\{\sum\limits_{k=0}^{j-1}|a_k(\tau)|\Big \slash \sum\limits_{k=j}^n a_k(\tau)\biggr\}\biggr\}, \phh M_{T,j}^*(t) \equiv \sist{M_{T,j}, \ph t \in [t_0,T],}{M_{t,j}, \ph t > T.}
$$

{\bf Lemma 2.3.} {\it Let for some $j=2,\ldots,n$ the inequalities $a_k(t) \ge 0, \ph k=\overline{j,n}, \linebreak \sum\limits_{k=j}^n a_k(t) >~ 0, \ph t \ge t_0$ be satisfied. Then $M_{T,j}^*(t)$ is a sub solution of the inequality (2.3) on $[t_0,+\infty)$.
}

Proof. It is obvious that $M_{T,j}^*(t)$ is a nondecreasing and continuous function on $[t_0,+\infty)$. Let $t_1 > t_0$ be fixed. To prove the lemma it is enough to show that $\eta_{t_1}(t) \equiv M_{t,j}, \ph t\in [t_0,t_1]$ is a solution of the inequality (2.3) on $[t_0,t_1]$. Since $a_k(t) \ge 0, \ph k=\overline{j,n}, \ph t \ge t_0$ we have $\sum\limits_{k=j}^n a_k(t)\eta^k \ge \Bigl[\sum\limits_{k=j}^n a_k(t)\Bigr]\eta^j$ for all $\eta \ge 1,$ and $t\ge t_0$. Then under the restriction $A_j(t) \equiv \sum\limits_{k=j}^n a_k(t) > 0, \ph t \ge t_0$ we get
$$
\sum\limits_{k=0}^n a_k(t) \eta^k \ge A_j(t)\eta^j\biggl[1 - \sum\limits_{k=0}^{j-1}|a_k(t)|\Big\slash A_j(t)\eta\biggr]
$$
for all $\eta \ge 1, \ph t \ge t_0.$ It follows from  here that $\eta_{t_1}(t)$ is solution of the inequality (2.3) on $[t_0,t_1]$. The lemma is proved.

For any $\gamma \in \mathbb{R}$ and $t_1\ge t_0$ we set
$$
\eta_{\gamma,t_1}(t) \equiv \gamma + \exp\biggl\{-\il{t_0}{t}a_1(\tau)d\tau\biggr\}\biggl[c(t_1) - \il{t_0}{t}\exp\biggl\{\il{t_0}{\tau}a_1(s) d s\biggr\} a_0(\tau) d\tau\biggr], \phh  t \in [t_0,t_1],
$$
where $c_1(t_1)\equiv \max\limits_{\xi\in[t_0,t_1]}\il{t_0}{\xi}\exp\biggl\{\il{t_0}{\tau}a_1(s) d s\biggr\} a_0(\tau) d\tau$.

{\bf Lemma 2.4.} {\it Let the following conditions be satisfied.

\noindent
(1) \ph $a_n(t) \ge 0, \ph t \ge t_0$.

\noindent
(2) \ph $a_k(t) = a_n(t) c_k(t) + d_k(t), \ph k=\overline{2,n-1}, \ph t \ge t_0$, where $c_k(t), \ph k=\overline{2,n-1}$ are bounded functions on $[t_0,t_1]$ for every $t_1 \ge t_0$ and

\noindent
(3) \ph $\sum\limits_{k=2}^{n-1} d_k(t) u^k \ge 0, \ph u \ge N_T, \ph t \ge t_0$, where
$$
N_{t_1}\equiv \max\Bigl\{1, \sup\limits_{t \in [t_0,t_1]}\sum\limits_{k=2}^n|c_k(t)|\Bigr\}, \ph t_1 \ge T, \ph\mbox{for some}  \ph T\ge t_0.
$$
Then
$$
\eta_T^*(t) \equiv \sist{\eta_{N_T,T}(t), \ph t\in [t_0,T],}{\eta_{N_{t},t}(t), \ph t \ge T}
$$
is a sub solution of the inequality (2.3) on $[t_0,+\infty)$.
}

Proof. Obviously, $\eta_T^*(t) \in C([t_0,+\infty))$. Therefore, to prove the lemma it is enough to show that for every $t_1 \ge T$ the function $\eta_{N_{t_1},t_1}(t)$ is a solution of the inequality (2.3) on $[t_0,t_1]$ and $\eta_{N_{t_1},t_1}(t_0) \ge \eta_T^*(t_0)$. The last inequality follows immediately from the definition of $\eta_{N_{t_1},t_1}(t)$. Consider the function
$$
F(t,u)\equiv 1 + \fr{c_{n-1}(t)}{u}+ \ldots + \fr{c_2(t)}{u^{n-2}}, \phh t \in [t_0,t_1], \ph u \ge 1 \phh (t_1 \ge T).
$$
Obviously $F(t,u) \ge 1 - \fr{\sum\limits_{k=2}^{n-1}|c_k(t)|}{u} \ge 0$ for all $t\in[t_0,t_1]$ and for all $u \ge N_{t_1}$. Moreover, $\eta_{N_{t_1},t_1}(t) \ge N_{t_1}, \phh t \in [t_0,t_1].$ Hence,
$$
F(t,\eta_{N_{t_1},t_1}(t)) \ge 0, \phh t \in [t_0,t_1]. \eqno (2.8)
$$
It is clear that $\eta_{N_{t_1},t_1}(t) \in C^1([t_0,t_1])$ and $\eta_{N_{t_1},t_1}'(t) + a_1(t)\eta_{N_{t_1},t_1}(t) + a_0(t) = 0, \ph  t \in [t_0,t_1].$ it follows from here and from the condition (2) that
$$
\eta_{N_{t_1},t_1}'(t) + \sum\limits_{k=0}^n a_k(t)\eta_{N_{t_1},t_1}^k(t) = a_n(t))\eta_{N_{t_1},t_1}^n(t) F(t,\eta_{N_{t_1},t_1}(t)) + \sum\limits_{k=2}^{n-1} d_k(t)\eta_{N_{t_1},t_1}^k(t),
$$
$t \in [t_0,t_1].$ This together with the conditions (1), (2) and with the inequality (2.8) implies that $\eta_{N_{t_1},t_1}(t)$ is a solution of the inequality (2.3) on $[t_0,t_1]$. The lemma is proved.

We set
$$
\eta_c(t)\equiv \exp\biggl\{-\il{t_0}{t}a_1(\tau) d\tau\biggr\}\biggl[c - \il{t_0}{t}\exp\biggl\{\il{t_0}{\tau}a_1(s) d s\biggr\}a_0(\tau) d\tau\biggr], \phh t \ge t_), \phh c \in \mathbb{R}.
$$

{\bf Lemma 2.5.} {\it Let the following conditions be satisfied.

\noindent
(4) \ph $a_2(t) > 0, \ph t \in [t_0,T],$

\noindent
(5) \ph  for some  $c \ge \max \limits_{t \in[t_0,T]}\il{t_0}{t}\exp\biggl\{\il{t_0}{\tau}a_1(s)d s\biggr\}a_0(\tau) d\tau$ the inequality \\  \phantom{aaaaaaaaaa}$\sum\limits_{k=3}^n|a_k(t)|\eta_c^{k-2}(t) \le~ a_2(t), \ph t \in [t_0,T]$ is valid.

\noindent
Then the function $\eta_c(t)$ is a nonnegative solution of the inequality (2.3) on $[t_0,T]$.
}

Proof. Obviously
$$
\eta_c(t) \ge 0, \phh t \in [t_0,T] \eqno (2.9)
$$
and
$$
\eta_c'(t) + a_1(t) \eta_c(t) + a_0(t) = 0, \phh t \in [t_0,T]. \eqno (2.10)
$$
It follows from the conditions (4), (5) and from the inequality (2.9) that
$$
\sum\limits_{k=2}^n a_k(t) \eta_c^k(t) = a_2(t) \eta_c^2(t)\Biggl[1 + \frac{\sum\limits_{k=3}^n a_k(t)\eta_c^{k-2}(t)}{ a_2(t)}\Biggr] \ge a_2(t) \eta_c^2(t)\Biggl[1 - \frac{\sum\limits_{k=3}^n a_k(t)\eta_c^{k-2}(t)}{ a_2(t)}\Biggr] \ge 0,
$$
$t\in [t_0,T]$. This together with (2.9) and (2.10) implies that $\eta_c(t)$ is a nonnegative solution of the inequality (2.3) on $[t_0,T]$. The lemma is proved.

We set
$$
\alpha(t) \equiv \sum\limits_{k=2}^n|a_k(t)| - a_1(t), \ph \theta_c(t) \equiv \exp\biggl\{{t_0}{t}\alpha(\tau) d\tau\biggr\}\biggl[c - \il{t_0}{t}\exp\biggl\{-\il{t_0}{\tau}\alpha(s) d s\biggr\} a_0(\tau) d\tau\biggr],
$$
$t \ge t_0, \ph c\in \mathbb{R}.$

{\bf Lemma 2.6.} {\it Let for some $c \ge \max\limits_{t \in [t_0,T]}\il{t_0}{t}\exp\biggl\{-\il{t_0}{\tau}\alpha(s) d s\biggr\}a_0(\tau) d\tau$ the inequality $\theta_c(t) \le 1, \ph t \in [t_0,T]$ be satisfied. Then $\theta_c(t)$ is a nonnegative solution of the inequality (2.3) on $[t_0,T].$
}

Proof. It is obvious that
$$
\theta_c(t) \ge 0, \phh t \in [t_0,T]. \eqno (2.11)
$$
Show that $\theta_c(t)$ satisfies (2.3) on $[t_0,T]$. We have
$$
\sum\limits_{k=0}^n a_k(t) \theta_c^k(t) = \Bigl(\sum\limits_{k=2}^n|a_k(t)|\Bigr)\theta_c(t) + \sum\limits_{k=2}^n a_k(t) \theta_c^k(t) + a_0(t) - \alpha(t) \theta_c(t), \ph t \in [t_0,T]. \eqno (2.12)
$$
Obviously,
$$
\theta_c'(t) + a_0(t) - \alpha(t) \theta_c(t) = 0, \phh t \in [t_0,T]. \eqno (2.13)
$$
It follows from here and from (2.12) that if $\sum\limits_{k=2}^{n}|a_k(t)| = 0$ for some fixed $t \in [t_0,T]$, then $\theta_c(t)$ satisfies (2.3) in $t$. Assume $\sum\limits_{k=2}^{n}|a_k(t)| \ne 0$ for a fixed $t \in [t_0,T]$. Then it follows from the condition $\theta_c(t) \le 1, \ph t \in [t_0,T]$ of the lemma and from (2.11) that
$$
\Bigl(\sum\limits_{k=2}^n|a_k(t)|\Bigr)\theta_c(t) + \sum\limits_{k=2}^n a_k(t) \theta_c^k(t) \ge
\Bigl(\sum\limits_{k=2}^n|a_k(t)|\Bigr)\theta_c(t)\Biggl[1 - \frac{\sum\limits_{k=2}^n|a_k(t)|\theta_c^2(t)}{\sum\limits_{k=2}^n|a_k(t)|\theta_c(t)}\Biggr] \ge 0
$$
for that fixed $t$. This together with (2.12) and (2.13) implies that $\theta_c(t)$ satisfies (2.3) in that fixed $t$. Hence, $\theta_c(t)$ satisfies (2.3) for all $t\in [t_0,T]$. The lemma is proved.

 Let $F(t,Y)$ be a continuous in $t$ and continuously differentiable in $Y$ vector function on $[t_0,\pinf)\times \mathbb{R}^m$. Consider the nonlinear system
$$
Y' = F(t,Y), \phh t \ge t_0. \no (2.14)
$$
Every solution $Y(t) = Y(t,t_0,Y_0)$ of this system exists either only a finite interval $[t_0,T)$ or is continuable on $[t_0,\pinf)$

\vsk

{\bf Lemma 2.7([19, p. 204, Lemma]).} {\it If a solution $Y(t)$ of the system (2.14) exists only on a finite interval $[t_0,T)$, then
$$
||Y(t)|| \to \pinf \ph \mbox{as} \ph t \to T- 0,
$$
where $||Y(t)||$ is any euclidian norm of $Y(t)$ for every fixed $t \in [t_0,T)$.
}

\phantom{aaaaaaaaaaaaaaaaaaaaaaaaaaaaaaaaaaaaaaaaaaaaaaaaaaaaaaaaaaaaaaaaaaaaaaa} $\blacksquare$

{\bf Lemma 2.8.} {\it For $k$ odd the inequality
$$
S_k(u,v) \ge 0, \ph u,v \in \mathbb{R} \ph \mbox{is valid}.
$$
}

Proof. If $u=0$, then
$$
S_k(u,v) = v^{k-1} = v^{2m} \ge 0, \phh v \in \mathbb{R}, \ph (m \in \mathbb{Z}_+). \eqno (2.15)
$$
For $u\ne 0$ we have
$$
S_k(u,v) = u^{2m} P_k(x), \phh x \equiv \frac{v}{u}, \phh P_k(x) \equiv \sum\limits_{j=0}^{k-1} x^j, \ph x \in \mathbb{R}. \eqno (2.16)
$$
Since $k-1$ is even all roots of $P_k(x)$ are complex (not real). Besides $P_k(9) =1 > 0$. Hence, $P_k(x) > 0, \ph x \in \mathbb{R}.$ This together with (2.15) and (2.16) implies that $S_k(u,v) \ge 0$ for all $u,v \in \mathbb{R}.$. The lemma is proved.

{\bf Lemma 2.9} {\it For $k$ even the inequality
$$
\frac{\partial S_k(u,v)}{\partial u} \ge 0, \phh u, v \in \mathbb{R} \phh\mbox{is valid}.
$$
}

Proof. Since $\frac{\partial S_k(u,v)}{\partial u} = (k-1) u^{k-2} + \ldots + 2 u v^{k-3} + v^{k-2}, \ph u,v \in R$ and $k$ is even, we have
$$
\frac{\partial S_k(u,v)}{\partial u}\Bigl|_{u=0} = (k-1) u^{2m}\ge 0, \phh u\in \mathbb{R}, \phh (m\in \mathbb{Z}_+). \eqno (2.17)
$$
For $u \ne 0$ the following equality is valid
$$
\frac{\partial S_k(u,v)}{\partial u} = u^{2m}Q_k(x), \phh x \equiv \frac{v}{u}, \phh Q_k(x) \equiv \sum\limits_{j=0}^{k-1}(j+1) x ^j, \phh x \in \mathbb{R}. \eqno (2.18)
$$
Consider the polynomials $q_j(x) = (j+1) x^{2j}(1 + x^2), \ph j=0,1, \ldots$. Obviously,
$$
q_j(x) \ge 0, \phh x \in \mathbb{R}, \phh j=0,1,2,\ldots .\eqno (2.19)
$$
and
$$
q_0(x) + (k-1) x^{2m} \ge 0, \phh x \in \mathbb{R}, \phh m \in \mathbb{Z}_+. \eqno (2.20)
$$
It is not difficult to verify that
$$
Q_k(x) = q_0(x) + \ldots + q_{2(m-1)}(x) +  (k-1) x^{2m}, \phh x \in \mathbb{R}, \phh m \in \mathbb{Z}_+.
$$
This together with (2.17)-(2.20) implies that $\frac{\partial S_k(u,v)}{\partial u} \ge 0, \ph u, v \in \mathbb{R}.$ The lemma is proved.

\vskip 10pt

Let $f(t,u)$ be a real-valued continuous function on $[t_0,T]\times \mathbb{R}$. Consider the first order differential equation
$$
y' = f(t,y), \phh t \in [t_0,T] \eqno (2.21)
$$
and the differential inequalities
$$
\zeta' \le f(t,\zeta), \phh t \in [t_0,T], \eqno (2.22)
$$
$$
\eta' \ge f(t,\eta), \phh t \in [t_0,T]. \eqno (2.23)
$$

\vskip 10pt

{\bf Theorem 2.1 ([20, Theorem 2.1])} {\it Let $\zeta(t)$ and $\eta(t)$ be solutions of the inequalities (2.22) and (2.23) respectively on $[t_0,T]$ such that $\zeta(t) \le \eta(t), \ph t\in [t_0,T], \ph  \zeta(t_0) \le~ \zeta(T), \linebreak \eta(t_0) \ge \eta(T)$. If any
solution $y(t)$ of the Cauchy problem
$y = f (t, u),  \ph y(t_0) = y_0 \in [\zeta(t_0),\eta(t_0)]$
 is unique,  then Eq. (2.1) has a solution $y_*(t)$ on $[t_0,T]$ such that $y_*(t_0) = ~y_*(T), \linebreak  \zeta(t) \le y_*(t) \le \eta(t), \ph t \in[t_0,T]$
}

Let  $ t_0 < t_1 < \ldots$ be a finite or infinite sequence such that   $t_k\in[t_0,\tau_0],\phantom{a}k=0,1,...$.

{\bf Definition 2.3.} {\it The sequence $\{t_k\}$
we will  call an usable sequence for the interval $[t_0,\tau_0)$, if   the maximum of the numbers $t_k$ coincides with  $\tau_0$ for finite $\{t_k\}$ , and  $\lim\limits_{k\to\infty} t_k = \tau_0$ for infinite $\{t_k\}$.
}

Let $a(t), \ph b(t)$ and $c(t)$ be real valued continuous functions on $[t_0,\tau_0) \ph (\tau \le +\infty)$. Consider the Riccati equation
$$
y' + a(t) y^2 + b(t) y + c(t) = 0, \ph t \in [t_0,\tau_0). \eqno (2.24)
$$

\vskip 10pt

{\bf Theorem 2.2 ([2, Theorem 4.1])}. {\it Assume  $a(t)\ge 0,\phantom{a} t\in [t_0,\tau_0)$  and
$$
\int\limits_{t_k}^t\exp\left\{\int\limits_{t_k}^\tau\biggl[b(s) - a(s)\biggl(\il{t_k}{s}\exp\biggl\{-\il{\xi}{s} b(\zeta) d\zeta\biggr\}c(\xi) d\xi\biggr)\biggr]ds\right\}
c(\tau)d\tau\le 0,\phantom{aa} t\in [t_k,t_{k+1}),
$$
$ k = 1,2, ...\phantom{a}$, where $\{t_k\}$ is an usable sequence for $[t_0,\tau_0)$.
Then for every  $\gamma \ge 0$ Eq.  (2.24) has a solution  $y_0(t)$ on  $[t_0,\tau_0)$, satisfying the initial condition  $y_0(t_0)=\gamma$, and $y_0(t)\ge 0,\phantom{a} t\in [t_0,\tau_0)$.
}

\vskip 10pt

{\bf Remark 2.1.} {\it Theorem 2.2 remains valid if for $\tau_0 < +\infty$ we replace $[t_0,\tau_0)$ by $[t_0,\tau_0]$ in it.}

We set
$$
I_\gamma(t) \equiv \gamma\exp\biggl\{-\il{t_0}{t} a_1(\tau) d\tau\biggr\} + \il{t_0}{t}\exp\biggl\{-\il{\tau}{t} a_1(s) d s\biggr\}|a_0(\tau)|d\tau, \phh t \ge t_0.
$$

\vskip 10pt

{\bf Theorem 2.3.} {\it Assume  $a_k(t)\ge 0,\phantom{a} k=\overline{2,n}, \ph  t\in [t_0,\tau_0)$  and
$$
\int\limits_{t_k}^t\exp\left\{\int\limits_{t_k}^\tau\biggl[a_1(s) - a_2(s)\biggl(\il{t_k}{s}\exp\biggl\{-\il{\xi}{s} a_1(\zeta) d\zeta\biggr\}a_0(\xi) d\xi\biggr)\biggr]ds\right\}
a_0(\tau)d\tau\le 0, \eqno (2.25)
$$
$\phantom{a} t\in [t_k,t_{k+1}), \ph   k = 1,2, ...$  where $\{t_k\}$ is an usable sequence for $[t_0,\tau_0)$..
Then for every  $\gamma \ge 0$ the inequality (2.3) has a solution  $\eta_\gamma^0(t)$ on  $[t_0,\tau_0)$, satisfying the initial condition  $\eta_\gamma^0(t_0)=\gamma$, and $0 \le \eta_\gamma^0(t) \le I_\gamma(t), \ph t\ge t_0$.
}

Proof. By Theorem 2.2 it follows from the conditions  $a_2(t)\ge 0, \ph  t\in [t_0,\tau_0)$ and (2.25) that for every $\gamma\ge 0$ every solution $y_\gamma(t)$ of the Riccati equation
$$
y' + a_2(t) y^2 + a_1(t) y + a_0(t) = 0, \phh t \in [t_0,\tau_0) \eqno (2.26)
$$
with $y_\gamma(t_0) = \gamma$ exists on $[t_0,\tau_0)$ and is nonnegative. It follows from here and from the condition
$a_k(t)\ge 0,\phantom{a} k=\overline{2,n}, \ph  t\in [t_0,\tau_0)$ of the theorem that $\eta_\gamma^0(t) \equiv y_\gamma(t)$ is a nonnegative solution of the inequality (2.3) on $[t_0,\tau_0)$ for every $\gamma \ge 0$.
Note that we can interpret $y_\gamma(t)$ as a solution of the linear equation
$$
x' + [a_2(t)y_\gamma(t) + a_1(t)]x + a_0(t) = 0, \ph t \ge t_0.
$$
Then by the Cauchy formula we have
$$
y_\gamma(t) = \gamma\exp\biggl\{-\il{t_0}{t}[a_2(\tau)y_\gamma(\tau) + a_1(\tau)]d\tau\biggr\} - \il{t_0}{t}\exp\biggl\{-\il{\tau}{t}[a_2(s) y_\gamma(s) + a_1(s)]d s\biggr\}a_0(\tau) d \tau,
$$
$t \ge t_0.$ Hence, $0\le \eta_\gamma^0(t) = y_\gamma(t) \le I_\gamma(t), \ph t \ge t_0.$ The theorem is proved.

\vskip 10pt

{\bf 3. Comparison criteria.} In this section we prove two comparison criteria for Eq. (1.1). These criteria with the aid of section 2 we use in section 4 to obtain some global solvability criteria for Eq. (1.1).

\vskip 10pt

{\bf Theorem 3.1.} {\it Let $y_1(t)$ be a solution of Eq. (2.1) on $[t_0,+\infty)$ and let $\eta^*(t)$ be a sub solution of the inequality (2.3) on $[t_0,+\infty)$ such that $y_1(t_0) < \eta^*(t_0)$. Moreover, let the following conditions be satisfied

\noindent
$(I) \ph D(t,u,y_1(t)) \le D_1(t,u,y_1(t)), \ph u \ge y_1(t), \ph t \ge t_0$, where $D_1(t,u,y_1(t))$ is a nondecreasing in $u \ge y_1(t)$ function for every $t \ge t_0$.

\noindent
$(II) \ph \gamma - y_1(t) + \il{t_0}{t}\exp\biggl\{\il{t_0}{\tau}D_1(s,\eta^*(s),y_1(s)) d s\biggr\}\Biggl(\sum\limits_{k=0}^n[b_k(\tau) - a_k(\tau)] y_1^k(\tau)\biggr) d\tau \ge 0, \phh t \ge t_0$ for some $\gamma\in[y_1(t_0),\eta^*(t_0)]$.

\noindent
Then  every solution $y(t)$ of Eq. (1.1) with $y(t_0) \in [\gamma,\eta^*(t_0)]$ exists on $[t_0,+\infty)$ and
$$
y_1(t) \le y(t) \le \eta^*(t), \phh t \ge t_0.
$$
Furthermore, if $y_1(t_0) < y(t_0) \ph (y(t_0) < \eta^*(t_0))$, then
$$
y_1(t) < y(t) \ph (y(t) < \eta^*(t)), \ph t \ge t_0.
$$
}

Proof. Let $y(t)$ be a solution of Eq. (1.1) with $y(t_0) \in [\gamma,\eta^*(t_0)]$ and let $[t_0,t_1)$ be its maximum existence interval. Then by Lemma 2.1 we have
$$
y(t) \le \eta^*(t), \ph t \in [t_0,t_1), \eqno (3.1)
$$
and if $y(t_0) < \eta^*(t_0),$ then
$$
y(t) < \eta^*(t), \phh t \in[t_),t_1). \eqno (3.2)
$$
In virtue of (2.2) we have
$$
y(t) - y_1(t) = \exp\biggl\{-\il{t_0}{t}D(\tau,y(\tau),y_1(\tau))d\tau\biggr\}\biggl[y(t_0) - y_1(t_0) - \phantom{aaaaaaaaaaaaaaaaaaaaaaaa}
$$
$$
-\il{t_0}{t}\exp\biggl\{\il{t_0}{\tau}D(s,y(s),y_1(s))d s\biggr\}\biggl(\sum\limits_{k=0}^n[a_k(\tau) - b_k(\tau)]y_1^k(\tau)\biggr)d\tau\biggr], \ph t \in[t_0,t_1). \eqno (3.3)
$$
Let us show that
$$
y_1(t) \le y(t), \ph t\in[t_0,t_1). \eqno (3.4)
$$
At first we consider the case $y(t_0) > y_1(t_0)$. Show that in this case
$$
y_1(t) < y(t), \phh t \in [t_0,t_1). \eqno (3.5)
$$
Suppose it is not true. Then there exists $t_2\in (t_0,t_1)$ such that
$$
y_1(t) < y(t), \phh t \in [t_0,t_2).
$$
$$
y_1(t_2) = y(t_2). \eqno (3.6)
$$
It follows from here, from (3.1) and from the condition $(I)$ that
$$
D(t,y(t),y_1(t)) \le D_1(t,\eta^*(t),y_1(t)), \phh t\in[t_0,t_2).
$$
Hence, the function
$$
H(\tau)\equiv \exp\biggl\{\il{t_0}{\tau}\Bigl[D(s,y(s),y_1(s))-D_1(s,\eta^*(s),y_1(s))\Bigr]ds\biggr\}, \phh \tau \in[t_0,t_2)
$$
is positive and non increasing on $[t_0,t_2)$. By mean value theorem for integrals (see [21, p. 869]) it follows from here that
$$
\il{t_0}{t}\exp\biggl\{\il{t_0}{\tau}D(s,y(s),y_1(s))d s\biggr\}\biggl(\sum\limits_{k=0}^n[a_k(\tau) - b_k(\tau)]y_1^k(\tau)\biggr) d\tau = \phantom{aaaaaaaaaaaaaaaaaaaaaaaaa}
$$
$$
\phantom{aaaaaaaaaaaaaaaa}= \il{t_0}{\kappa(t)}\exp\biggl\{\il{t_0}{\tau}D_1(s,\eta^*(s),y_1(s))d s\biggr\}\biggl(\sum\limits_{k=0}^n[a_k(\tau) - b_k(\tau)]y_1^k(\tau)\biggr) d\tau
$$
for some $\kappa(t) \in [t_0,t), \ph t \in [t_0,t_2)$. This together with (3.3) and with the condition $(II)$ implies that $y_1(t_2) < y(t_2)$, which contradicts (3.6). The obtained contradiction proves (3.5), hence proves (3.4). Let us show that (3.4) is also valid in the case $y(t_0) = y_1(t_0).$ Suppose, for some $t_3 \in(t_0,t_1)$
$$
y(t_3) < y_1(t_3). \eqno (3.7)
$$
Let $\widetilde{y}_\delta(t)$ be a solution of Eq. (1.1) with  $\widetilde{y}_\delta(t_0) > y_1(t_0)$. Then by already proven (3.5) we have $\widetilde{y}_\delta(t_3) > y_1(t_3)$. As far as the solutions of Eq. (1.1) continuously depend on their initial values we chose $\delta > 0$ so small that  $\widetilde{y}_\delta(t_3) - y_1(t_3) < \frac{y_1(t_3) - y(t_3)}{2}$. Since, $\widetilde{y}_\delta(t_3) > y_1(t_3)$ it follows from the last inequality that
$$
y_1(t_0) - y(t_3) < \widetilde{y}_\delta(t_3) - y_1(t_3) < \frac{y_1(t_3) - y(t_3)}{2},
$$
which contradicts (3.7). The obtained contradiction proves that (3.4) is also valid for $y(t_0) = y_1(t_0).$ Note that the proof of (3.3) and (3.4) in the general case $y(t_0) \ge y_1(t_0)$ repeats the proof of them for the case $y(t_)) > y_1(t_0)$. Therefore, due to (3.1), (3.2), (3.4) and (3.5) to complete the proof of the theorem it remains to show that
$$
t_1 = +\infty. \eqno (3.8)
$$
Suppose $t_1 , +\infty$. Then it follows from (3.1) and (3.4) that $y(t)$ is bounded on $[t_0,t_1)$. By lemma 2.7 it follows from here that $[t_0,t_1)$ is not the maximum existence interval for $y(t)$, which contradicts our supposition. The obtained contradiction proves (3.8). The proof of the theorem is completed.

Note that every function $y_1(t) \equiv \zeta(t) \in C^1([t_0,+\infty))$ is a solution of Eq. (2.1) with $b_0(t) = -\zeta'(t), \ph b_1(t)=\ldots=b_n(t)\equiv 0, \ph t \ge t_0$. Then
$$
\sum\limits_{k=0}^n[b_k(t) - a_k(t)]y_1^k(t) = - \Bigl[\zeta'(t) + \sum\limits_{k=0}^na_k(t)\zeta^k(t)\Bigr] \phh t \ge t_0.
$$
From here and from Theorem 3.1 we obtain immediately

\vskip 10pt

{\bf Corollary 3.1.} {\it Let $\eta^*(t)$ be a sub solution of the inequality (2.3) on $[t_0,+\infty)$ and let for some $y_1(t) \equiv \zeta(t) \in C^1([t_0,+\infty))$ with $\zeta(t_0) < \eta^*(t_0)$ the condition $(I)$ of Theorem~ 3.1 and the following condition be satisfied
$$
(II^0) \ph \zeta(t_0) - \gamma + \il{t_0}{t}\exp\biggl\{\il{t_0}{\tau}D_1(s,\eta^*(s),\zeta(s)) d s\biggr\}\biggl(\zeta'(\tau) + \sum\limits_{k=0}^n a_k(\tau)\zeta^k(\tau)\biggr)d\tau \le 0, \ph t \ge t_0,
$$
for some $\gamma \in [\zeta(t_0),\eta^*(t_0)].$

\noindent
Then every solution $y(t)$ of Eq. (1.1) with $y(t_0) \in [\gamma,\eta^*(t_0)]$ exists on $[t_0,+\infty)$ and
$$
\zeta(t) \le y(t) \le \eta^*(t), \phh t \ge tt_0.
$$
Furthermore, if $\zeta(t_0) < y(t_0)\ph (y(t_0) < \eta^*(t_0)),$ then
$$
\zeta(t) < y(t) \ph (y(t) < \eta^*(t)), \phh t \ge t_0.
$$
}

\phantom{aaaaaaaaaaaaaaaaaaaaaaaaaaaaaaaaaaaaaaaaaaaaaaaaaaaaaaaaaaaaaaa} $\blacksquare$

\vsk
{\bf Remark 3.1.} {\it It is clear form the proofs of Theorem 3.1 and Corollary 3.1 that we can replace $\eta^*(t)$ in the conditions $(II)$ and $(II^0)$ respectively  of Theorem 3.1 and Corollary~ 3.1 by a continuous function $\widetilde{\eta}^*(t) \ge \eta^*(t), \ph t \in [t_0,+\infty)$.
}

Let $e_k(t), \ph k=\overline{0,n}$ be real-valued continuous functions on $[t_0,+\infty)$. Consider the equation
$$
y' + \sum\limits_{k=0}^n e_k(t) y^k = 0, \phh t \ge t_0. \eqno (3.9)
$$

\vsk

{\bf Theorem 3.2} {\it Let $y_1(t)$ and $y_2(t)$ be solutions of the equations (2.1) and (3.9) respectively on $[t_0,+\infty)$ such that $y_1(t_0)\le y_2(t_0)$ and let the following conditions be satisfied.

\noindent
$(III) \ph \sum\limits_{k=0}^n(b_k(t) - a_k(t)) y_1^k(t) \ge 0, \ph t \ge t_0$,

\noindent
$(IV) \ph  \sum\limits_{k=0}^n(e_k(t) - a_k(t)) y_2^k(t) \le 0, \ph t \ge t_0$.

\noindent
Then every solution $y(t)$  of Eq. (1.1) with $y(t_0)\in [y_1(t_0),y_2(t_0)]$ exists on $[t_0,+\infty)$ and
$$
y_1(t) \le y(t) \le y_2(t), \phh t \ge t_0.
$$
Furthermore, if $y_1(t_0) < y(t_0) \ph (y(t_0) < y_2(t_0)),$ then
$$
y_1(t) < y(t) \phh (y(t) < y_2(t)), \phh t \ge t_0.
$$
}

Proof. Let $y(t)$ be a solution of Eq. (1.1) with $y(t_0) \in [y_1(t_0),y_2(t_0)]$ and let $[t_0,t_1)$ be its maximum existence interval. Then by (2.2) the following equations are valid
$$
y(t) - y_1(t) = \exp\biggl\{-\il{t_0}{t}D(\tau,y(\tau),y_1(\tau))d\tau\biggr\}\biggl[y(t_0) - y_1(t_0) - \phantom{aaaaaaaaaaaaaaaaaaaaaaa}
$$
$$
\phantom{aaaaaaa}-\il{t_0}{t}\exp\biggl\{\il{t_0}{\tau}D(s,y(s),y_1(s))ds\biggr\}\biggl(\sum\limits_{k=0}^n[a_k(\tau) - b_k(\tau)]y_1^k(\tau)\biggr)d\tau, \phh t \in[t_0,t_1),
$$
$$
y(t) - y_2(t) = \exp\biggl\{-\il{t_0}{t}D(\tau,y(\tau),y_2(\tau))d\tau\biggr\}\biggl[y(t_0) - y_2(t_0) - \phantom{aaaaaaaaaaaaaaaaaaaaaaa}
$$
$$
\phantom{aaaaaaa}-\il{t_0}{t}\exp\biggl\{\il{t_0}{\tau}D(s,y(s),y_2(s))ds\biggr\}\biggl(\sum\limits_{k=0}^n[a_k(\tau) - e_k(\tau)]y_1^k(\tau)\biggr)d\tau, \phh t \in[t_0,t_1).
$$
It follows from here and from the conditions $(III)$ and $(IV)$ of the theorem that
$$
y_1(t) \le y(t) \le y_2(t), \phh t \in [t_0,t_1). \eqno (3.10)
$$
and if $y_1(t_0) < y(t_0) \ph (y(t_0) < y_2(t_0))$, then
$$
y_1(t) < y(t) \ph (y(t) < y_2(t)), \phh t \in [t_0,t_1).
$$
Therefore, to complete the proof of the theorem it remains to show that
$$
t_1 = +\infty. \eqno (3.11)
$$
Suppose $t_1<+\infty$. Then by Lemma 2.7 it follows from (3.10) that $[t_0,t_1)$ is not the maximum existence interval for $y(t)$, which contradicts our supposition. The obtained contradiction proves (3.11). The proof of the theorem is completed.

\vskip 10pt

{\bf Corollary 3.2.} {\it Let $\eta^*(t)$ and $\zeta^*(t)$ be sub and super solutions of the inequalities (2.3) and (2.4) respectively on $[t_0,+\infty)$ such that $\zeta^*(t) \le \eta^*(t)$. Then every solution $y(t)$ of Eq. (1.1) with $y(t_0) \in [\zeta^*(t_0),\eta^*(t_0)]$ exists on $[t_0,+\infty)$ and
$$
\zeta^*(t) \le y(t) \le \eta^*(t), \phh t \ge t_0.
$$
Furthermore, if $\zeta^*(t_0) < y(t_0) \ph (y(t_0) < \eta^*(t_0))$, then
$$
\zeta^*(t) < y(t) \phh (y(t) < \eta^*(t)), \phh t \ge t_0.
$$
}

Proof. To prove the corollary it is enough to show that for every $\tau_0> t_0$ and solutions $\zeta(t)$ and $\eta(t)$ of the inequalities (2.4) and (2.3) respectively on $[t_0,\tau_9]$ with $\zeta(t_0)\le \eta(t_0)$ any solution $y(t)$ of Eq. (1.1) with $y(t_0) \in [\zeta(t_0),\eta(t_0)]$ exists on $[t_0,\tau_0]$ and
$$
\zeta(t) \le y(t) \le \eta(t), \phh t \in [t_0,\tau_0], \eqno (3.12)
$$
and if $\zeta(t_0) < y(t_0) \ph (y(t_0) < \eta(t_0))$, that
$$
\zeta(t) < y(t) \phh (y(t) < \eta(t)), \phh [t_0,\tau_0]. \eqno (3.13)
$$
The function $y_1(t)\equiv \zeta(t)$ is a solution of Eq. (2.1) on $[t_0,\tau_0]$ for $b_0(t)\equiv -\zeta'(t), \ph b_1(t)=,\ldots = b_n(t) \equiv 0$, and $y_2(t) \equiv \eta(t)$ is a solution of the equation (3.9) on $[t_0,\tau_0]$ for $e_0(t)\equiv -\eta'(t), \ph e_1(t)=\ldots = e _n(t) \equiv 0$. Then the condition $(III)$ gives us
$$
\zeta'(t) + \sum\limits_{k=0}^n a_k(t) \zeta^k(t) \le 0, \phh t \ge t_0
$$
and the condition $(IV)$ gives us
$$
\eta'(t) + \sum\limits_{k=0}^n a_k(t) \eta^k(t) \ge 0, \phh t \ge t_0.
$$
Therefore by Theorem 3.2 (note that Theorem 3.2  remains valid if we replace $[t_0,+\infty)$ by $[t_0,\tau_0]$ in it) then the inequalities (3.12) and (3.13) are valid. The corollary is proved.

\vsk

{\bf 4. Global solvability criteria}. Denote $a_k^+(t)\equiv \max\{0,a_k(t)\}, \ph k=\overline{2,n}, \ph t \ge t_0$.

\vskip 10pt

{\bf Theorem 4.1.} {\it Let the condition of Lemma 2.3 and the  following condition be satisfied

\noindent
$(A)$ \ph for a nonnegative $\zeta(t) \in C^1([t_0,+\infty))$ with $\zeta(t_0) < M_{\gamma,T}^*(t_0)$ and for some \linebreak $\nu\in [\zeta(t_0),M_{\gamma,T}^*(t_0)]$
$$
\zeta(t_0) - \nu +\il{t_0}{t}\exp\biggl\{\il{t_0}{\tau}\Bigl[\sum\limits_{k=2}^na_k^+(s) S_k(M_{\gamma,T}^*(s),\zeta(s)) + a_1(s)\Bigr]d s\biggr\}\times \phantom{aaaaaaaaaaaaaaaaaaaaaaaaaaaaaaaaaa}
$$
$$
\phantom{aaaaaaaaaaaaaaaaaaaaaaaaaAAAAAaaaaaa}\times\biggl(\zeta'(\tau) + \sum\limits_{k=0}^n a_k(\tau)\zeta^k(\tau)\biggr)d \tau \le 0, \phh t \ge t_0.
$$

\noindent
Then every solution $y(t)$ of Eq. (1.1) with $y(t_0) \in [\nu,M_{\gamma,T}^*(t_0)]$ exists on $[t_),+\infty)$ and
$$
\zeta(t) \le y(t) \le M_{\gamma,T}^*(t), \phh t \ge t_0. \eqno (4.1)
$$
Furthermore, if $\zeta(t_0) < y(t_0) \ph (y(t_0) < M_{\gamma,T}^*(t_0))$, then
$$
\zeta(t) < y(t) \ph (y(t) < M_{\gamma,T}^*(t)), \phh t \ge t_0. \eqno (4.2)
$$
}

Proof. By Lemma 2.1 $M_{\gamma,T}^*(t)$ is a sub solution of the inequality (2.3) on $[t_0,+\infty)$. Note that $y_1(t) \equiv \zeta(t)$ is a solution of Eq. (2.1) on $[t_0,+\infty)$ for $b_0(t) \equiv -\zeta'(t), \ph b_1(t)=\ldots = b_n(t) \equiv 0$. Then since $\zeta(t)$ is nonnegative we have
$$
D(t,u,v) \le \sum\limits_{k=2}^na_k^+(t)S_k(u,\zeta(t)) + a_1(t), \ph \mbox{for all} \ph u \ge \zeta(t), \ph t \ge t_0.
$$
Moreover, $\sum\limits_{k=2}^na_k^+(t)S_k(u,\zeta(t))$ is nondecreasing in $u\ge \zeta(t) \ge 0,$ for all $t \ge t_0$. It follows from here and from $(A)$ that the conditions of Corollary 3.1 are satisfied. Hence, every solution $y(t)$ of Eq. (1.1) with $y(t_0) \in [\nu,M_{\gamma,T}^*(t_0)]$ exists on $[t_0,+\infty)$ and the inequalities (4.1) and (4.2) are valid. The theorem is proved.

By analogy with the proof of Theorem 4.1 it can be proved the following theorem

\vskip 10pt
{\bf Theorem 4.2.} {\it Let the conditions of Lemma 2.4 and the condition

\noindent
$(B)$ for a nonnegative $\zeta(t) \in C^1([t_0,+\infty))$ with $\zeta(t_0) < \eta_T^*(t_0)$ and for some $\nu\in[\zeta(t_0),\eta_T^*(t_0)]$
$$
\zeta(t_0) -\nu + \il{t_0}{t}\exp\biggl\{\il{t_0}{\tau}\Bigl[\sum\limits_{k=2}^n a_k^+(s) S_k(\eta_T^*(s),\zeta(s)) + a_1(s)\Bigr]d s\biggr\}\times \phantom{aaaaaaaaaaaaaaaaaaaaaaaaaaaaaaaaaaa}
$$
$$
\phantom{aaaaaaaaaaaaaaaaaaaaaaaaaaaaaaaaaaaaaa}\times\biggl(\zeta'(\tau) + \sum\limits_{k=0}^n a_k(\tau)\zeta^k(\tau)\biggr)d\tau \le 0, \phh t\ge t_0.
$$
be satisfied.

\noindent
Then every solution $y(t)$ of Eq. (1.1) with $y(t_0) \in [\nu,\eta_T^*(t_0)]$ exists on $[t_0,+\infty)$ and
$$
\zeta(t) \le y(t) \le \eta_T^*(t), \phh t \ge t_0.
$$
Furthermore, if $\zeta(t_0) < y(t_0) \ph (y(t_0) < \eta_T^*(t_0))$, then
$$
\zeta(t) < y(t) \ph (y(t) < \eta_T^*(t)), \phh t\ge t_0.
$$
}

\phantom{aaaaaaaaaaaaaaaaaaaaaaaaaaaaaaaaaaaaaaaaaaaaaaaaaaaaaaaaaaaaaaaaaaaa} $\blacksquare$

\vsk

{\bf Corollary 4.1.} {\it Let the conditions of Lemma 2.3 or Lemma 2.4 be satisfied. If $a_0(t) \le~ 0, \ph t \ge t_0$, then every solution $y(t)$ of Eq. (1.1) with $y(t_0) \ge 0$ exists on $[t_0,+\infty)$ and is nonnegative
}

Proof. Let $y(t)$ be a solution of Eq. (1.1) with $y(t_0) \ge 0$. Under the conditions of Lemma 2.3 (of Lemma 2.4) we can take $M_T^*(t) \ph (\eta_T^*(t))$ so that $y(t_0) \le M_T^*(t) \ph (y(t_0) \le \eta_T^*(t))$. Then the condition
$a_0(t) \le~ 0, \ph t \ge t_0$ provides the satisfiability of the condition $(A)$ of Theorem 4.1 (of the condition $(B)$ of Theorem 4.2) for $\nu=0, \ph \zeta(t) \equiv 0$. Hence, the assertion of the corollary is valid. The corollary is proved.

Let $\{t_k\}$ be the same sequence for $\tau_0 = +\infty$, defined in the section 2.

\vsk
{\bf Theorem 4.3.} {\it Let the following conditions be satisfied.

\noindent
$ (C) \ph
a_k(t) \ge 0, \ph k=\overline{2,n},  \ph t \ge t_0.
$

\noindent
$ (D) \ph
\int\limits_{t_l}^t\exp\left\{\int\limits_{t_l}^\tau\biggl[a_1(s) - a_2(s)\biggl(\il{t_l}{s}\exp\biggl\{-\il{\xi}{s} a_1(\zeta) d\zeta\biggr\}a_0(\xi) d\xi\biggr)\biggr]ds\right\}
a_0(\tau)d\tau\le 0, \\
\phantom{aaaa}  t\in [t_l,t_{l+1}), \ph   l = 1,2, \ldots$, where $\{t_l\}$ is an usable sequence for $[t_0,+\infty)$..

\noindent
$ (E) \ph
\il{t_0}{t}\exp\biggl\{\il{t_0}{\tau}\Bigl[\sum\limits_{k=2}^n a_k(s)I_\gamma^{k-1}(s) + a_1(s)\Bigr]d s\biggr\} a_0(\tau)d\tau \le 0, \phh t\ge t_0.
$

\noindent
Then every solution $y(t)$ of Eq. (1.1) with $y(t_0) = \gamma \ge 0$ exists on $[t_0,+\infty)$ and
$$
0\le y(t) \le I_\gamma(t), \ph t \ge t_0.
$$
}

Proof. By Theorem 2.3 it follows from the conditions $(C)$ and $(D)$ of the theorem that for every $\gamma > 0$ the inequality (2.3) has a nonnegative solution $\eta_\gamma^0(t)$ with $\eta_\gamma^0(t_0) = \gamma$. It follows from the condition $(C)$ of the theorem that $D(t,u,0) \le D_1(t,u,0) \equiv \sum\limits_{k=1}^n a_k(t) u^{k-1}, \ph u \ge~ 0$ and $D_1(t,u,0)$ is a nondecreasing function for $u\ge 0$. Taking into account Remark 2.1 by Corollary 3.1 from here and from  $(E)$ we obtain that every solution $y(t)$ of Eq. (1.1) with $y(t_0)= \gamma \ge 0$ exists on $[t_0,+\infty)$ and $0\le y(t) \le I_\gamma(t), \ph t \ge t_0$. The theorem is proved.

\vsk

We set $\sigma_k^\pm \equiv \frac{1 \pm (-1)^k}{2}, \ph k=0,1,2,\ldots.$ Obviously,
$$
\sigma_k^+ = \sist{1, \ph \mbox{for} \ph  k \ph \mbox{even},}{0, \ph \mbox{for} \ph  k \ph \mbox{odd},} \phh \phh \sigma_k^- = \sist{0, \ph \mbox{for} \ph  k \ph \mbox{even},}{1, \ph \mbox{for} \ph  k \ph \mbox{odd},} \phh k=0,1,\ldots.
$$

\vskip 10pt

{\bf Theorem  4.4.} {\it Let the conditions of Lemma 2.4 and the following conditions be satisfied.

\noindent
$(F) \ph (-1)^k a_k(t) \ge 0, \ph k = \overline{2,n}, \phh t \ge t_0,$

\noindent
$(G)$ \ph for some $\zeta(t) \in C^1([t_0,+\infty))$ with $\zeta(t_0) < \eta_T^*(t_0)$ and for some $\nu \in [\zeta(t_0),\eta_T^*(t_0)]$
$$
\zeta(t_0) - \nu + \il{t_0}{t}\exp\biggl\{\il{t_0}{\tau}\Bigl[\sum\limits_{k=2}^n\sigma_k^+ a_k^+(s) S_k(\eta_T^*(s),\zeta(s)) + a_1(s)\Bigr]d s\biggr\} \times \phantom{aaaaaaaaaaaaaaaaaaaaaaaaaaaaaaaaaaa}
$$
$$
\phantom{aaaaaaaaaaaaaaaaaaaaaaaaaaaaaa}\times\biggl(\zeta'(\tau) + \sum\limits_{k=0}^n a_k(\tau)\zeta^k(\tau)\biggr) d \tau \le 0, \phh t \ge t_0.
$$

\noindent
Then every solution $y(t)$ of Eq. (1.1) with $y(t_0) \in [\nu,\eta_T^*(t_0)]$ exists on $[t_0,+\infty)$ and
$$
\zeta(t) \le y(t) \le \eta_T^*(t), \phh t \ge t_0. \eqno (4.3)
$$
Furthermore, if $\zeta(t_0) < y(t_0) \ph (y(t_0) < \eta_T^*(t_0))$, then
$$
\zeta(t) < y(t) \ph (y(t) < \eta_T^*(t)), \phh t \ge t_0. \eqno (4.4)
$$
}

Proof. By virtue of Lemma 2.4 $\eta_T^*(t)$ is a sub solution of the inequality (2.3) on $[t_0,+\infty)$. Since
$$
D(t,u,v) = \sum\limits_{k=2}^n\sigma_k^+ a_k(t) S_k(u,v) + \sum\limits_{k=2}^n\sigma_k^- a_k(t) S_k(u,v) + a_1(t), \phh u, v \in \mathbb{R}, \ph t \ge t_0
$$
By Lemmas 2.8 and 2.9 it follows from $(F)$ that $\sum\limits_{k=2}^n\sigma_k^+ a_k(t) S_k(u,v) + a_1(t)$ is nondecreasing in $u\ge \zeta(t)$ for all $t$ and $\sum\limits_{k=2}^n\sigma_k^- a_k(t) S_k(u,v) \le 0, \ph t \ge t_0$.
Hence, $D(t,u,\zeta(t)) \le \sum\limits_{k=2}^n \sigma_k^+ a_k(t) S_k(u,\zeta(t)) + a_1(t), \ph u \ge \zeta(t), \ph t \ge t_0.$ It follows from here and from $(G)$ that the condition $(II)$ of Theorem 3.1 is satisfied for the case $b_0(t) = -\zeta'(t), \ph b_1(t) = \ldots = b_n(t) \equiv 0$. Thus, all conditions of Theorem 3.1 are satisfied. Therefore, every solution $y(t)$ of Eq. (1.1) with $y(t_0) \in [\nu,\eta_T^*(t_0)]$ exists on $[t_0,+\infty)$ and the inequalities (4.3) and (4.4) are satisfied. The theorem is proved.

{\bf Theorem 4.5.} {\it Let the conditions of Lemma 2.5 and the following  condition be satisfied.

\noindent
$(H)$ \ph for a nonnegative $\zeta(t) \in C^1([t_0,T])$ with $\zeta(t_0) < \eta_c(t_0)$ and for sone $\nu \in [\zeta(t_0),\eta_c(t)]$
$$
\zeta(t_0) - \nu + \il{t_0}{t}\exp\biggl\{\il{t_0}{\tau}\Bigl[\sum\limits_{k=2}^n a_k^+(s) S_k(\eta_c(s),\zeta(s)) + a_1(s)\Bigr] d s\biggr\}\times \phantom{aaaaaaaaaaaaaaaaaaaaaaaaa}
$$
$$
\phantom{aaaaaaaaaaaaaaaaaaaaaaaaaaaaaa}\times\biggl(\zeta'(\tau) + \sum\limits_{k=0}^n a_k(\tau) \zeta^k(\tau)\biggr) d \tau  \le 0, \phh t \in [t_0,T].
$$
Then every solution $y(t)$ of Eq. (1.1) with $y(t_0) \in [\nu,\eta_c(t_0)]$ exists on $[t_0,T]$ and
$$
\zeta(t) \le y(t) \le \eta_c(t), \phh t \in [t_0,T]. \eqno (4.5)
$$
Furthermore, if $\zeta(t_0) < y(t_0) \ph (y(t_0) < \eta_c(t_0))$, then
$$
\zeta(t) < y(t) \ph (y(t) < \eta_c(t)), \phh t \in [t_0,T]. \eqno (4.6)
$$
}

Proof. By virtue of Lemma 2.5 $\eta_c(t)$ is a solution of the inequality (2.3) on $[t_0,T]$. Since $\zeta(t)$ is nonnegative we have
$$
D(t,u,\zeta(t)) \le \sum\limits_{k=2}^n a_k^+(t) S_k(u,\zeta(t)) + a_1(t), \phh t \in [t_0,T].
$$
It follows from here and from the condition $(H)$ that the condition $(II)$ of Theorem 3.1 for $[t_0,T]$ and for the case $b_0(t)\equiv -\zeta'(t), \ph b_1(t) = \ldots =b_n(t) \equiv 0, \ph t \in [t_0,T]$ is satisfied. Thus all conditions of Theorem 3.1 for $[t_0,T]$ are satisfied. Therefore, every solution $y(t)$ of Eq. (1.1) with $y(t_0) \in [\nu,\eta_c(t_0)]$ exists on $[t_0,T]$ and the inequalities (4.5) and (4.6) are satisfied. The theorem is proved.

By analogy with the proof of Theorem 4.5 it can be proved the following theorem

{\bf Theorem 4.6.} {\it Let the condition of Lemma 2.6 and the following condition be satisfied

\noindent
 for a nonnegative $\zeta(t) \in C^1([t_0,T])$ with $\zeta(t_0) < \theta_c(t_0)$ and for some $\nu \in [\zeta(t_0),\theta_c(t_0)]$
$$
\zeta(t_0) - \nu + \il{t_0}{t}\exp\biggl\{\il{t_0}{\tau}\Bigl[\sum\limits_{k=2}^n a_k^+(s) S_k(\theta_c(s),\zeta(s)) + a_1(s)\Bigr]d s\biggr\}\times \phantom{aaaaaaaaaaaaaaaaaaaaaaaaaaaaaaaaaa}
$$
$$
\phantom{aaaaaaaaaaaaaaaaaaaaaaaaaaaaaaaaa} \times\biggl(\zeta'(\tau) + \sum\limits_{k=0}^n a_k(\tau) \zeta^k(\tau)\biggr)d\tau \le 0, \ph t \in [t_0,T].
$$

\noindent
Then every solution $y(t)$ of Eq. (1.1) with $y(t_0) \in [\nu,\theta_c(t_0)]$ exists on $[t_0,T]$ and
$$
\zeta(t) \le y(t) \le \theta_c(t), \phh t \in [t_0,T].
$$
Furthermore, if $\zeta(t_0) < y(t_0), \ph (y(t_0) < \theta_c(t_0))$, then
$$
\zeta(t) < y(t), \ph (y(t) < \theta_c(t)), \phh t \in [t_0,T].
$$
}

\phantom{aaaaaaaaaaaaaaaaaaaaaaaaaaaaaaaaaaaaaaaaaaaaaaaaaaaaaaaaaaaaaaaaaaaa} $\blacksquare$

{\bf Corollary 4.2.} {\it Let the conditions of Lemma 2.5 and the following conditions be satisfied

\noindent
$(I) \ph a_1(t) < 0, \ph t\in [t_0,T]$,

\noindent
$(J)$ \ph for some $\zeta_0  \in (0,\eta_c(t_0))$ with $\sum\limits_{k=2}^n|a_k(t)|\zeta_0^{k-1} \le |a_1(t)|, \ph t \in [t_0,T]$ and  for some $\nu \in [\zeta_0,\eta_c(t_0)]$
$$
\zeta_0 - \nu + \il{t_0}{t}\exp\biggl\{\il{t_0}{\tau}\Bigl[\sum\limits_{k=2}^n a_k^+(s) S_k(\eta_c(s),\zeta_0) + a_1(s)\Bigr] d s\biggr\}a_0(\tau) d \tau \le 0, \phh t \in [t_0,T].
$$

\noindent
Then every solution $y(t)$ of Eq. (1.1) with $y(t_0) \in [\nu,\eta_c(t_0)]$ exists on $[t_0,T]$ and
$$
\zeta_0 \le y(t) \le \eta_c(t), \phh t \in [t_0,T]. \eqno (4.7)
$$
Furthermore, if $\zeta_0 < y(t_0), \ph (y(t_0) < \eta_c(t_0))$, then
$$
\zeta_0 < y(t), \ph (y(t) < \eta_c(t)), \phh t \in [t_0,T]. \eqno (4.8)
$$
}

Proof. It follows from the condition $(I)$ that for some (enough small) $\zeta_0 \in (0,\eta_c(t_0))$ with $\sum\limits_{k=2}^n|a_k(t)|\zeta_0^{k-1} \le |a_1(t)|, \ph t \in [t_0,T]$   the inequality $\sum\limits_{k=2}^n a_k(t) \zeta_0^k \le 0, \ph t \in [t_0,T]$ is satisfied. This together with the condition $(J)$ implies the condition $(H)$ of Theorem 4.5. Thus all conditions of Theorem 4.5 are satisfied. Therefore, every solution $y(t)$ of Eq. (1.1) with $y(t_0) \in [\nu,\eta_c(t_0)]$ exists on $[t_0,T]$ and the inequalities (4.7) and (4.8) are satisfied. The corollary is proved.

By analogy with the proof of Corollary 4.2 one can prove the following assertion.

\vsk

{\bf Corollary 4.3.} {\it Let the conditions of Lemma 2.5 and the following conditions be satisfied

\noindent
$ \ph a_1(t) > 0, \ph t\in [t_0,T]$,

\noindent
 for some $\zeta_0  < 0$ with $\sum\limits_{k=2}^n|a_k(t)||\zeta_0|^{k-1} \le a_1(t), \ph t \in [t_0,T]$ and  for some $\nu \in [\zeta_0,\eta_c(t_0)]$
$$
\zeta_0 - \nu + \il{t_0}{t}\exp\biggl\{\il{t_0}{\tau}\Bigl[\sum\limits_{k=2}^n a_k^+(s) S_k(\eta_c(s),\zeta_0) + a_1(s)\Bigr] d s\biggr\}a_0(\tau) d \tau \le 0, \phh t \in [t_0,T].
$$

\noindent
Then every solution $y(t)$ of Eq. (1.1) with $y(t_0) \in [\nu,\eta_c(t_0)]$ exists on $[t_0,T]$ and
$$
\zeta_0 \le y(t) \le \eta_c(t), \phh t \in [t_0,T].
$$
Furthermore, if $\zeta_0 < y(t_0), \ph (y(t_0) < \eta_c(t_0))$, then
$$
\zeta_0 < y(t), \ph (y(t) < \eta_c(t)), \phh t \in [t_0,T].
$$
}

\phantom{aaaaaaaaaaaaaaaaaaaaaaaaaaaaaaaaaaaaaaaaaaaaaaaaaaaaaaaaaaaaaaaaaa} $\blacksquare$

For any $\gamma \in \mathbb{R}, \ph t_1 \ge t_0$ we set
$$
\zeta_{\gamma,t_1}(t) \equiv -\gamma -\exp\biggl\{-\il{t_0}{t}a_1(\tau) d \tau\biggr\}\biggl[c(t_1) + \il{t_0}{t}\exp\biggl\{\il{t_0}{\tau} a_1(s) d s\biggr\}a_0(\tau) d \tau\biggr], \ph t \in [t_0,t_1),
$$
where $c(t_1) \equiv \max\limits_{\xi \in[t_0,t_1]}\biggl(-\il{t_0}{\xi}\exp\biggl\{\il{t_0}{\tau}a_1(s) d s\biggr\}a_0(\tau) d \tau\biggr).$

\vsk

{\bf Theorem 4.7.} {\it Let the conditions of Lemma 2.4 and the following conditions be satisfied.

\noindent
$(K) \ph \sum\limits_{k=2}^{n-1}(-1)^{k+1} d_k(t) u^k \ge 0,$ for all $u\ge N_T$ and $t \ge t_0$.

\noindent
$(L) \ph n$ is odd.

\noindent
Then every solution $y(t)$ of Eq. (1.1) with $y(t_0) \in [\zeta_T^*(t_0),\eta_T^*(t_0)]$ exists on $[t_0,+\infty)$ and
$$
\zeta_T^*(t) \le y(t) \le \eta_T^*(t), \phh t \ge t_0, \eqno (4.9)
$$
where $\eta_T^*(t)$ is defined in Lemma 2.4 and
$\zeta_T^*(t) \equiv \sist{\zeta_{N_T,T}(t), \ph t \in [t_0,T],}{\zeta_{N_t,t}(t), \ph t \ge T,}$ such that $\zeta_T^*(t_0) \le \eta_T^*(t_0).$ furthermore, if $\zeta_T^*(t_0) < y(t_0) \ph (y(t_0) < \eta_T^*(t_0))$, then
$$
\zeta_T^*(t) < y(t) \ph (y(t) < \eta_T^*(t)), \phh t \ge t_0. \eqno (4.10)
$$
}

Proof. By Lemma 2.4 $\eta_T^*(t)$ is a sub solution of the inequality (2.3) on $[t_0,+\infty)$. Show that $\zeta_T^*(t)$ is a super solution of the inequality (2.4) on $[t_0,+\infty)$. Consider the differential inequality
$$
\eta' + \sum\limits_{k=0}^n \widetilde{a}_k(t)\eta^k \ge 0, \phh t \ge t_0, \eqno (4.11)
$$
where $\widetilde{a}_k(t) = (-1)^{k+1} a_k(t), \ph k=\overline{0,n}, \ph t \ge t_0$ It follows from $(K)$ and from the condition $(1)$ of Lemma 2.4 that

\noindent
$\widetilde{(1)} \ph \widetilde{a}_n(t) \ge 0, \ph t \ge t_0.$

\noindent
it follows from the condition $(2)$ of Lemma 2.4 that

\noindent
$\widetilde{(2)} \ph \widetilde{a}_k(t) = \widetilde{a}_n(t) \widetilde{c}_k(t) + \widetilde{d}_k(t), \ph k=\overline{2,n-1}, \ph t \ge t_0$, where $\widetilde{c}_k(t) = (-1)^{k+1} c_k(t), \ph k=\overline{2,n-1}, \ph t \ge t_0$ are bounded function on $[t_0,t_1]$ for every $t_1 \ge t_0, \ph \widetilde{d}_k(t) = (-1)^{k+1} d_k(t), \ph k=\overline{2,n-1}, \ph t \ge t_0$.

\noindent
It follows from the condition $(K)$, that

\noindent
$\widetilde{(3)} \ph \sum\limits_{k=2}^{n-1}\widetilde{d}_k(t) u^k \ge 0$ for all $u\ge N_T, \ph t \ge t_0$.

\noindent
We see that all conditions of Lemma 2.4 for the inequality (4.11) are satisfied. Hence, by Lemma 2.4 $\widetilde{\eta}_T^*(t) \equiv -\zeta_T^*(t)$ is a sub solution of the inequality (4.11) on $[t_0,+\infty)$. Then $\zeta_T^*(t)$ is a super solution of the inequality (2.4) on $[t_0,+\infty)$. By Corollary 3.2 it follows from here that every solution $y(t)$ of Eq. (1.1) with $y(t_0) \in [\zeta_T^*(t_0), \eta_T^*(t_0)]$ \ph (note that always  $\zeta_T^*(t_0)\le  \eta_T^*(t_0)]$) exists on $[t_0,+\infty)$ and the inequalities (4.9) and (4.10) are satisfied. The theorem is proved.

We set
$$
\theta_c^-(t) \equiv -\exp\biggl\{-\il{t_0}{t}\alpha(\tau) d\tau\biggr\}\biggl[c + \il{t_0}{t}\exp\biggl\{-\il{t_0}{\tau}\alpha(s) d s\biggr\} a_0(\tau) d \tau\biggr], \phh t \ge t_0, \phh t \in \mathbb{R}.
$$

\vsk
{\bf Theorem 4.8.} {\it Assume for some $c^+ \ge \max\limits_{t \in [t_0,T]}\il{t_0}{t}\exp\biggl\{-\il{t_0}{\tau}\alpha(s) d s\biggr\} a_0(\tau) d\tau, \linebreak c^- \ge - \max\limits_{t \in [t_0,T]}\il{t_0}{t}\exp\biggl\{-\il{t_0}{\tau}\alpha(s) d s\biggr\} a_0(\tau) d\tau$ the inequalities
$$
\theta_{c^+}(t) \le 1, \phh |\theta_{c^-}^-(t)| \le 1, \phh t \in [t_0,T]
$$
are valid. Then every solution $y(t)$ of Eq. (1.1) with $y(t_0) \in [\theta_{c^-}^-(t),\theta_{c^+}(t)]$ exists on $[t_0,T]$ and
$$
\theta_{c^-}^-(t) \le y(t) \le \theta_{c^+}(t), \phh t \in [t_0,T]. \eqno (4.12)
$$
Furthermore, if $\theta_{c^-}^-(t_0) < y(t_0) \ph (y(t_0) < \theta_{c^+}(t_0))$, then
$$
\theta_{c^-}^-(t) < y(t) \ph (y(t) < \theta_{c^+}(t)), \phh t \in [t_0,T]. \eqno (4.13)
$$
}

Proof. We have

\noindent
$
\theta_{c^-}^-(t_0) \le \min\limits_{t \in [t_0,T]}\il{t_0}{t}\exp\biggl\{-\il{t_0}{\tau}\alpha(s) d s\biggr\} a_0(\tau) d\tau \le \max\limits_{t \in [t_0,T]}\il{t_0}{t}\exp\biggl\{-\il{t_0}{\tau}\alpha(s) d s\biggr\} a_0(\tau) d\tau \le \theta_{c^+}(t_0).
$

\noindent
Therefore, the relation $y(t_0) \in [\theta_{c^-}^-(t),\theta_{c^+}(t)]$ is correct.  By Lemma 2.6 $\theta_{c^+}(t)$ is a solution of the inequality (2.3) on $[t_0,T]$, and $-\theta_{c^-}^-(t)$ is a solution of the inequality (4.11) on $[t_0,T]$. Then, since $\theta_{c^-}^-(t)\le \theta_{c^+}(t)]$, by Corollary 3.3 every solution $y(t)$ of Eq. (1.1) with $y(t_0) \in [\theta_{c^-}^-(t),\theta_{c^+}(t)]$ exists on $[t_0,T]$ and the inequalities (4.12) and (4.13) are valid. The theorem is proved.
\vsk

{\bf 5. \ph Closed solutions.}

\vsk

{\bf Theorem 5.1.} {\it Let the following conditions be satisfied.

\noindent
$1^0)$ \ph for some $j=2,\ldots,n$ the inequalities $a_k(t) \ge 0, \ph k=\overline{j,n}, \ph \sum\limits_{k=j}^n a_k(t)> 0, \ph t \in [t_0,T]$ are valid,

\noindent
$2^0)$ \ph $\il{t_0}{t}\exp\biggl\{\il{t_0}{\tau}\Bigl[\sum\limits_{k=2}^na_k^+(s)M_{T+\gamma}^{k-1} + a_1(s)\Bigr]d s\biggr\} a_0(\tau) d\tau \le 0, \ph t\in [t_0,T]$ \ph for some $\gamma \ge 0$.

\noindent
Then the following statements are valid.

\noindent
$\alpha) \ph$ Eq. (1.1) has a nonnegative closed solution $y_*(t)$  on $[t_0,T]$,

\noindent
$\beta)$ \ph If, in particular,  $a_0(t) \not \equiv 0$ and $a_0(t) \ge 0, \ph t \in [t_0,T],$
  then  $y_*(t)$ is positive,

\noindent
$\gamma)$ \ph If, in particular, $j=2$ and $\il{t_0}{T} a_1(\tau) d\tau > 0$,  then  $y_*(t)$ is isolated.
}

Proof. Let us prove $\alpha)$. It follows from the conditions of the theorem that for $\zeta(t) \equiv 0$ the conditions of Theorem 4.1 are satisfied. Then by Theorem 4.1 the solutions $y_1(t)$ and $y_2(t)$ of Eq. (1.1) with $y_1(t_0) =0$ and $y_2(t_0) = M_{T+\gamma}$ exist on $[t_0,T]$ and $y_1(T) \ge 0, \ph y_2(T) \le M_{\gamma,T}^*(T) = M_{\gamma,T}^*(t_0)$. By  Theorem 2.1 it follows from here that Eq. (1.1) has a nonnegative closed solution $y_*(t)$  on $[t_0,T]$. The statement $\alpha)$  is proved. Let us prove $\beta)$.
If  $a_0(t) \not \equiv 0$ and $a_0(t) \ge 0, \ph t \in [t_0,T],$
    then $y_1(t)\not \equiv 0$. Hence,
$$
y_1(t_1)> 0 \ph \mbox{for sone} \ph t_1\in [t_0,T] \eqno (5.1)
$$
(as far as $y_1(t) \ge 0, \ph t \in [t_0,T]$).
Consider the equation
$$
y' +\sum\limits_{k=1}^n a_k(t) y^k =0, \ph t \in [t_0,T].
$$
Since $y_0(t)\equiv 0$ is a solution of this equation by (2.2) we have
$$
y_1(T) = \exp\biggl\{-\il{t_1}{T}D(\tau,0,y_1(\tau)) d\tau\biggr\}\biggl[y_1(t_1) - \il{t_1}{T}\exp\biggl\{\il{t_1}{\tau}D(s,0,y_1(s))d s\biggr\} a_0(\tau) d \tau\biggr].
$$
It follows form here,  from  (5.1)   and from the conditions of $\beta)$ that $y_1(T) > 0 = y_1(t_0)$. Therefore, $y_1(t)$ is not a closed solution of Eq. (1.1) on $[t_0,T]$. By the uniqueness theorem it follows from here and from the statement $\alpha)$ that $y_*(t)$ is positive . The statement $\beta)$ is proved. It remains to prove $\gamma)$. Let us show that $y_*(t)$ is  isolated.   Suppose  $y_*(t)$ is not isolated.   Then there exists a sequence $\{y_m(t)\}_{m+1}^\infty$ of closed solutions of Eq. (1.1) on $[t_0,T]$ such that $y_m(t_0) \to y_*(t_0)$ for $m \to \infty$.  By (2.2) we have
$$
y_*(T) - y_m(T) = \exp\biggl\{-\il{t_0}{T}\Bigl(\sum\limits_{k=2}^n a_k(\tau) S_k(y_*(\tau),y_m(\tau)) + a_1(\tau)\Bigr) d\tau\biggr\}[y_*(t_0) - y_m(t_0)], \eqno (5.2)
$$
$m=1,2,\ldots.$ Since $j=2$ it follows from the conditions $1^0)$ that \linebreak $\il{t_0}{T}\Bigl(\sum\limits_{k=2}^n a_k(\tau) S_k(y_*(\tau),y_*(\tau))\Bigr) d \tau \ge~ 0$. Then as far as the solutions of Eq. (1.1) continuously depend on their initial values and $\il{t_0}{T}a_1(\tau) d \tau > 0$ we can chose $m=m_0$ so large that
$
\il{t_0}{T}\Bigl(\sum\limits_{k=2}^n a_k(\tau) S_k(y_*(\tau),y_m(\tau))+ a_1(\tau)\Bigr) d\tau > 0.
$
It follows from here  and from (5.2) that $y_{m_0}(t)$ is not closed. We have obtained a contradiction, proving $\gamma)$. The proof of the theorem is completed.

\vsk

{\bf Corollary 5.1.} {\it Let for some $j=2,\ldots,n$ the inequalities $(-1)^k a_k(t) \ge 0, \ph k=\overline{j,n}, \ph \sum\limits_{k=j}^n(-1)^k a_k(t) > 0, \ph a_0(t) \le 0, \ph t \in [t_0,T]$ be satisfied.

\noindent
Then the following statements are valid

\noindent
$\alpha^0)$ Eq. (1.1) has a non positive closed solution $y_*(t)$ on $[t_0,T]$.

\noindent
$\beta^0)$ \ph If, in particular,  $a_0(t)\not \equiv 0$ and $a_0(t) \le 0, \ph t \in [t_0,T]$, then  $y_*(t)$ is negative,

\noindent
$\gamma^0)$ \ph If, in particular, $j=2$ and $\il{t_0}{T} a_1(\tau) d\tau > 0$,  then $y_*(t)$ is isolated.

}

Proof. In Eq. (1.1) we substitute
$$
y = - z, \phh t \to - t. \eqno (5.3)
$$
We obtain
$$
z' + \sum\limits_{k=0}^n (-1)^k a_k(-t) z^k =0, \ph t \le -t_0.
$$
Then by Theorem 5.1 it follows from the conditions of the corollary that the transformed (last) equation has a  nonnegative closed solution $z_*(t)$ on $[-T,-t_0]$, for which the statements   $\alpha) - \gamma)$ of Theorem 5.1 are valid. It follows from here and from (5.3) that $y_*(t) \equiv -z_*(-t)$ is a nonnegative closed solution of Eq. (1.1) on $[t_0,T]$, for which the statements $\alpha^0) - \gamma^0)$ are valid. The corollary is proved.

Note that despite of Theorem 1.1 in the statement $\alpha)$ of  Corollary 5.1 the condition $a_0(t) \equiv 0$ of Theorem 1.1 is weakened up to $a_0(t) \le 0, \ph t \in [t_0,T]$ and the condition $\il{t_0}{T} a_1(\tau) d\tau > 0$ is omitted. Therefore, Corollary 5.1 is a complement of Theorem 1.1.

The inequality $\sum\limits_{k=j}^{n}a_k(t) > 0, \ph t\in [t_0,T]$ in conditions of  Theorem 5.1  looks like a strict limitation. The next theorem attempts to partially weaken it.

\vsk

{\bf Theorem 5.2.} {\it Let the conditions of Theorem 4.3 be satisfied. If $\sum\limits_{k=2}^na_k(t)\not \equiv 0$ or
$\il{t_0}{T}a_1(\tau) d \tau > 0$, then Eq. (1.1) has a nonnegative closed solution on $[t_0,T]$. In the case $\il{t_0}{T}a_1(\tau) d \tau > 0$ it is isolated.
}

Proof. By Theorem 4.3 for every $\gamma \ge 0$ Eq. (1.1) has a nonnegative solution $y_\gamma(t)$ on $[t_0,T]$ with $y_\gamma(t_0) = \gamma.$ Let us show that there exists $\gamma > 0$ such that
$$
y_\gamma(t_0) \ge y_\gamma(T). \eqno (5.4)
$$
At first we show that if $\sum\limits_{k=2}^na_k(t) \not \equiv 0, \ph t \in [t_0,T]$, then
$$
\lim\limits_{\gamma \to +\infty}\il{t_0}{T}\Bigl[\sum\limits_{k=2}^na_k(t)y_\gamma^{k-1}(t)\Bigr] d t = +\infty. \eqno (5.5)
$$
By (1.1) we can interpret $y_\gamma(t)$ as a solution of the linear equation
$$
x' +\Bigl[\sum\limits_{k=1}^na_k(t)y_\gamma^{k-1}(t)\Bigr] x + a_0(t) = 0, \phh t \in [t_0,T].
$$
Then by the Cauchy formula we have
$$
y_\gamma(t) = \gamma \exp\biggl\{-\il{t_0}{t}\Bigl[\sum\limits_{k=1}^na_k(\tau)y_\gamma^{k-1}(\tau)\Bigr]d\tau\biggr\} - \phantom{aaaaaaaaaaaaaaaaaaaaaaaaaaaaaaaaaa}
$$
$$
\phantom{aaaaaaaaaaaa}-\il{t_0}{t}\exp\biggl\{-\il{\tau}{t}\Bigl[\sum\limits_{k=1}^na_k(s)y_\gamma^{k-1}(s)\Bigr]d s\biggr\} a_0(\tau) d \tau,  \ph t \in [t_0,T]. \eqno (5.6)
$$
Multiplying both sides of this equality by $\Bigl[\sum\limits_{k=2}^na_k(t)y_\gamma^{k-2}(t)\Bigr]
\exp\biggl\{-\il{t_0}{t}\Bigl[\sum\limits_{k=2}^na_k(\tau)y_\gamma^{k-1}(\tau)\Bigr]d\tau\biggr\}$ and integrating over $[t_0,T]$ we obtain
$$
\exp\biggl\{\il{t_0}{T}\Bigl[\sum\limits_{k=2}^na_k(\tau)y_\gamma^{k-1}(\tau)\Bigr]d\tau\biggr\} = 1 + \gamma \il{t_0}{T}\Bigl[\sum\limits_{k=2}^na_k(t)y_\gamma^{k-2}(t)\Bigr]\exp\biggl\{-\il{t_0}{t}a_1(\tau) d \tau\biggr\} - \phantom{aaaaaaaaaa}
$$
$$
\phantom{aaaaaaaa}-\il{t_0}{T}\Bigl[\sum\limits_{k=2}^na_k(t)y_\gamma^{k-2}(t)\Bigr]d t \il{t_0}{t}\exp\biggl\{\il{t_0}{\tau}\Bigl[\sum\limits_{k=2}^na_k(s)y_\gamma^{k-1}(s)\Bigr]d s - \il{\tau}{t}a_1(s) d s\biggr\} a_0(\tau) d \tau
$$
From here we obtain
$$
\exp\biggl\{\il{t_0}{T}\Bigl[\sum\limits_{k=2}^na_k(\tau)y_\gamma^{k-1}(\tau)\Bigr]d\tau\biggr\} \ge 1 +
\il{t_0}{T}\Bigl[\sum\limits_{k=2}^na_k(t)y_\gamma^{k-2}(t)\Bigr]\exp\biggl\{-\il{t_0}{t}a_1(\tau) d \tau\biggr\} d t \times \phantom{aaaaaaa}
$$
$$
\phantom{aaaaaaaaaaa}\times
\biggl[\gamma -\il{t_0}{T} \exp\biggl\{\il{t_0}{\tau}\Bigl[\sum\limits_{k=2}^na_k(s)y_\gamma^{k-1}(s)\Bigr]d s + \il{t_0}{\tau}a_1(s) d s\biggr\} |a_0(\tau)| d \tau\biggr]. \eqno (5.7)
$$
Suppose
$$
\il{t_0}{T}\Bigl[\sum\limits_{k=2}^na_k(t)y_\gamma^{k-2}(t)\Bigr]d t  \le M, \ph \gamma > 0. \eqno (5.8)
$$
Then (5.7) implies
$$
\exp\biggl\{\il{t_0}{T}\Bigl[\sum\limits_{k=2}^na_k(\tau)y_\gamma^{k-1}(\tau)\Bigr]d\tau\biggr\} \ge 1 +
\il{t_0}{T}\Bigl[\sum\limits_{k=2}^na_k(t)y_\gamma^{k-2}(t)\Bigr]\exp\biggl\{-\il{t_0}{t}a_1(\tau) d \tau\biggr\} d t \times \phantom{aaaaaaa}
$$
$$
\phantom{aaaaaaaaaaaaaaaaaaaaaaaaaa}\times
\biggl[\gamma - \exp\biggl\{M\biggr\}\il{t_0}{T} \exp\biggl\{ \il{t_0}{\tau}a_1(s) d s\biggr\} |a_0(\tau)| d \tau\biggr]. \eqno (5.9)
$$
By the uniqueness theorem $y_\gamma(t) > 0, \ph t \in [t_0,T]$ for all $\gamma > 0$. Therefore (since $a_k(t) \ge~ 0, \linebreak t \in [t_0,T]$ and $\sum\limits_{k=2}^{n}a_k(t) \not \equiv 0$)

$$
\il{t_0}{T}\Bigl[\sum\limits_{k=2}^na_k(t)y_\gamma^{k-2}(t)\Bigr]\exp\biggl\{-\il{t_0}{t}a_1(\tau) d \tau\biggr\} d t \ge \il{t_0}{T}\Bigl[\sum\limits_{k=2}^na_k(t)y_{\gamma_0}^{k-2}(t)\Bigr]\exp\biggl\{-\il{t_0}{t}a_1(\tau) d \tau\biggr\} d t > 0
$$
for all $\gamma \ge \gamma_0 >0$. It follows from here that the right part of the inequality (5.9) tends to $+\infty$ as $\gamma \to +\infty$, whereas, according to (5.8) its lift part is bounded. We have obtained a contradiction, proving (5.5).  It follows from (5.6) that
$$
y_\gamma(T) = \gamma \exp\biggl\{-\il{t_0}{T}\Bigl[\sum\limits_{k=1}^na_k(\tau)y_\gamma^{k-1}(\tau)\Bigr]d\tau\biggr\} - \phantom{aaaaaaaaaaaaaaaaaaaaaaaaaaaaaaaaaa}
$$
$$
\phantom{aaaaaaaaaaaa}-\il{t_0}{T}\exp\biggl\{-\il{\tau}{t}\Bigl[\sum\limits_{k=1}^na_k(s)y_\gamma^{k-1}(s)\Bigr]d s\biggr\} a_0(\tau) d \tau,  \ph t \in [t_0,T], \ph \gamma \ge 0.
$$
Therefore $y_\gamma(T) \le y_\gamma(t_0) = \gamma$, provided
$$
\gamma\biggl(1 - \exp\biggl\{\il{t_0}{T}\biggl(\sum\limits_{k=1}^na_k(\tau) y_\gamma^{k-1}(\tau)\biggr)d\tau\biggr\}\biggr) \ge -  \il{t_0}{T}\exp\biggl\{-\il{\tau}{t}\Bigl[\sum\limits_{k=1}^na_k(s)y_\gamma^{k-1}(s)\Bigr]d s\biggr\} a_0(\tau) d \tau,
$$
which will be fulfilled if by virtue of (5.5) we chose $\gamma \ge  2 \il{t_0}{T}\exp\biggl\{\il{t_0}{\tau} a_1(s) d s\biggr\}|a_0(\tau)|d\tau$ so large that $\il{t_0}{T}\sum\limits_{k=1}^na_k(\tau)y_\gamma(\tau) d \tau \ge \ln 2$. Therefore (5.4) is  proved for the case  $\sum\limits_{k=2}^na_k(t)\not \equiv 0$. If $\sum\limits_{k=2}^na_k(t)\equiv 0$ and $\il{t_0}{T} a_1(\tau)d\tau > 0$, then from the obvious equality  $y_\gamma(T) = \gamma\exp\biggl\{-\il{t_0}{T}a_(\tau) d\tau\biggr\} -\il{t_0}{T}\exp\biggl\{\il{t_0}{\tau}a_1(s) d s\biggr\}a_0(\tau) d \tau$ we derive that for
$$
\gamma \ge \il{t_0}{T}\exp\biggl\{\il{t_0}{\tau}a_1(s) d s\biggr\}|a_0(\tau)| d \tau \bigg\slash \biggl(1 - \exp\biggl\{-\il{t_0}{T}a_1(\tau) d\tau\biggr\}\biggr)
$$
the inequality (5.4) is fulfilled. Thus, under the restriction   $\sum\limits_{k=2}^na_k(t)\not \equiv 0$ or $\il{t_0}{T} a_1(\tau)d\tau > 0$  of the theorem the inequality (5.4) is valid.
Then since $y_0(t_0) = 0 \le y_0(T),$ by Theorem~ 2.1 Eq. (1.1) has a nonnegative closed solution $y_*(t)$ on $[t_0,T]$.
To complete the proof of the theorem it remains to show that if  $\il{t_0}{T} a_1(\tau)d\tau > 0$, then $y_*(t)$ is isolated. The proof of this fact is similar to the proof of the assertion $\gamma)$ of Theorem 5.1. Therefore we omit it. The proof of the theorem is completed.

\vsk

{\bf Theorem 5.3.} {\it Let the following conditions be satisfied.

\noindent
$3^0) \ph a_n(t) \ge 0, \ph t \in [t_0,T]$,

\noindent
$4^0) \ph a_k(t) = a_n(t) c_k(t) + d_k(t), \ph k=\overline{2,n-1}$, where $c_k(t), \ph k=\overline{2,n-1}$ are bounded functions on $[t_0,T]$,

\noindent
$5^0) \ph \sum \limits_{k=2}^{n-1} d_k(t) u^k \ge 0, \ph u \ge N_T, \ph t \in [t_0,T]$,

\noindent
$6^0) \ph \max\limits_{t \in [t_0,T]}  \il{t_0}{t} \exp\biggl\{\il{t_0}{\tau}a_1(s) d s\biggr\} a_0(\tau) d\tau\biggl[1 - \exp\biggl\{\il{t_0}{T}a_1(\tau) d \tau\biggr\}\biggr] \le \il{t_0}{T}\exp\biggl\{\il{t_0}{\tau}a_1(s) d s \biggr\} a_0(\tau) d \tau.$

\noindent
$7^0) \ph \il{t_0}{t}\exp\biggl\{\il{t_0}{\tau}\Bigl[\sum\limits_{k=2}^n a_k^+(s)\eta_{\gamma,T}^{k-1}(s) + a_1(s)\Bigr]d s\biggr\} a_0(\tau) d\tau \le 0, \ph t \in [t_0,T].$ for some $\gamma\ge 0$,

\noindent

\noindent
Then Eq. (1.1) has a nonnegative closed solution on $[t_0,T]$.
}

Proof.
By virtue of Lemma 2.4 it follows from the conditions $3^0) - 5^0)$ that $\eta_{\gamma,T}(t)$ is a solution of the inequality (2.3) on $[t_0,T]$. Then it follows from the condition $7^0)$ that the conditions of Theorem 4.2 with $\zeta(t) \equiv 0$ are satisfied. Hence, according to Theorem~ 4.2 the solutions $y_1(t)$ and $y_2(t)$ of Eq. (1.1) with $y_1(t_0) = 0, \ph y_2(t_0) = \eta_{\gamma,T}(t_0) = c(T)$ exist on $[t_0,T]$ and $y_1(T) \ge 0, \ph y_2(T) \le \eta_{\gamma,T}(T)$. It follows from the condition $6^0)$ that $\eta_{\gamma,T}(T) \le \eta_{\gamma,T}(t_0)$. Therefore $y_2(T) \le y_2(t_0)$. By  Theorem 2.1 it follows from here that Eq. (1.1) has a nonnegative closed solution on $[t_0,T]$. The theorem is proved.

\vsk

{\bf Example 5.1.} {\it Consider the equation
$$
y' + \sum\limits_{k=0}^6 a_k(t) y^k = 0, \phh t \ge t_0, \eqno (5.10)
$$
where, $a_0(t) \equiv - \sin 10 t, \ph a_1(t)$ is any continuous function, $a_2(t) = \cos^4 t, \ph  a_3(t) = -2|sin t \cos^3 t|, \ph    a_4(t) = \sin^2 t \cos^2 t, \ph a_5(t) = - \sin^2 t |cos \pi t|, \ph  a_6(t) = \sin^2 t, \ph[t_0,T]$. Obviously,  the conditions of Corollary 5.1 for Eq. (5.10) are satisfied. It is not difficult to verify that the conditions of Theorem 5.3 with $c_2(t) =c_3(t)=c_4(t)\equiv0, \ph c_5(t) = -|\cos \pi t|, \ph d_2(t) = \cos^4 t, \ph d_3(t) = -2|\sin t \cos^3t|, \ph d(4) = \sin^2 t \cos^2 t, \ph d_5(t) \equiv 0, \ph t \in [t_0,T]$ for Eq. (5.10) are satisfied. Therefore Eq. (5.10) has at least a nonnegative closed solution $y_+(t)$ on $[t_0,T]$ and at least a non positive closed solution $y_-(t)$ on $[t_0,T]$ (for every $T > t_0)$. Since $a_0(t) \not\equiv 0$ we have $y_+(t) \ne y_-(t), \ph t \in [t_0,T]$.
}

\vsk

{\bf Theorem 5.4.} {\it Let the following conditions be satisfied.

\noindent
$8^0) \ph a_2(t) > 0, \ph t \in[t_0,T]$,

\noindent
$9^0)$ \ph for some $c \ge \max\limits_{t \in[t_0,T]}\il{t_0}{t}\exp\biggl\{\il{t_0}{\tau} a_1(s) d s\biggr\} a_0(\tau) d\tau$, the inequality $\sum\limits_{k=3}^n|a_k(t)|\eta_c^{k-2}(t) \le a_2(t), \ph t \in [t_0,T]$ is valid and

\noindent
$10^0) \ph \il{t_0}{t}\exp\biggl\{\il{t_0}{\tau}\Bigl[\sum\limits_{k=2}^na_k^+(s) \eta_c^{k-1}(s) + a_1(s)\Bigr]d s\biggr\} a_0(\tau) d \tau \le 0, \ph t \in [t_0,T]$,

\noindent
$11^0) \ph c\biggl(1 - \exp\biggl\{\il{t_0}{T}a_1(\tau) d\tau\biggr\}\biggr) \le \il{t_0}{T}\exp\biggl\{\il{t_0}{\tau}a_1(s) d s\biggr\} a_0(\tau) d\tau.$

\noindent
Then Eq. (1.1) has a nonnegative closed solution on $[t_0,T]$.
}

Proof. By Lemma 2.5 it follows from the conditions $8^0)$ and $9^0)$ that $\eta_c(t)$ is a solution of the inequality
(2.3) on $[t_0,T]$. It follows from the condition $10^0)$ that the condition $(E)$ with $\zeta(t)\equiv 0$ of Theorem 4.5 is satisfied. It follows from the condition $11^0)$ that $\eta_c(t_0) \ge \eta_c(T)$. Then by Theorems 2.1 and 4.5 Eq. (1.1) has a nonnegative closed solution on $[t_0,T]$. The theorem is proved.

Let us write $a_2(t) = \lambda p(t), \ph p(t) > 0, \ph t\in [t_0,T]$. Then for all   $\lambda \ge \lambda_0 \equiv \linebreak \equiv \max\limits_{t \in[t_o,T]}\{(\sum\limits_{k=3}^{n}|a_k(t)|\eta_c^{k-2}(t))/p(t)\}$ the condition $9^0)$ of Theorem 5.4  will be  satisfied. If we write $a_2(t) = \lambda + p(t), \ph p(t) \in C([t_0,T])$, then for all $\lambda \ge \lambda_0 \equiv \max\limits_{t \in[t_o,T]}\{(\sum\limits_{k=3}^{n}|a_k(t)|\eta_c^{k-2}(t)) - p(t)\}$ the condition $9^0)$ of the Theorem 5.4  will be  satisfied as well. Despite of this in Theorems 2 and 3 of work [18] the parameter $\lambda_0$ is undetermined. Moreover, for $a_0(t) \equiv 0, \ph c = 0$ the conditions $10^0)$ and $11^0)$ of Theorem 5.4  are satisfied. Therefore, Theorem 5.4 is a complement of both mentioned above Theorems 2 and 3.

\vsk
{\bf Theorem 5.5.} {\it Let the following conditions be satisfied.

\noindent
$12^0) \ph a_n(t) \ge 0, \ph t \in [t_0,T]$,

\noindent
$13^0) \ph a_k(t) = a_n(t) c_k(t) + d_k(t), \ph k=\overline{2,n-1}, \ph t \in [t_0,T],$ where $c_k(t), \ph k=\overline{2,n-1}$ are bounded functions on $[t_0,T]$ and

\noindent
$14^0) \ph \sum\limits_{k=2}^{n-1} d_k(t) u^k \ge 0, \ph u \ge N_T, \ph t\in [t_0,T],$

\noindent
$15^0) \ph \sum\limits_{k=2}^{n-1} (-1)^{k+1} d_k(t) u^k \ge 0,  \ph u \ge N_T, \ph  t\in [t_0,T],$

\noindent
$16^0) \ph n$ is odd,

\noindent
$17^0) \ph  \max\limits_{\xi \in[t_0,T]}\biggl(\il{t_0}{\xi}\exp\biggl\{\il{t_0}{\tau}a_1(s) d s\biggr\}a_0(\tau) d\tau\biggr) \biggl[1 - \exp\biggl\{\il{t_0}{T}a_1(\tau) d\tau\biggr\}\biggr] \le\\
\phantom{aaaaaaaaaaaaaaaaaaaaaaaaaaaaaaaaaaaaaaaaaaaaaaaaa}\le  \il{t_0}{T}\exp\bigg\{\il{t_0}{\tau}a_1(s) d s\biggr\}a_0(\tau) d\tau$,

\noindent
$18^0) \ph \min\limits_{\xi \in[t_0,T]}\biggl(\il{t_0}{\xi}\exp\biggl\{\il{t_0}{\tau}a_1(s) d s\biggr\}a_0(\tau) d\tau\biggr)\biggl[1 - \exp\biggl\{\il{t_0}{T}a_1(\tau) d\tau\biggr\}\biggr] \ge \\ \phantom{aaaaaaaaaaaaaaaaaaaaaaaaaaaaaaaaaaaaaaaaaaaaaaaaa} \ge \il{t_0}{T}\exp\bigg\{\il{t_0}{\tau}a_1(s) d s\biggr\}a_0(\tau) d\tau$.

\noindent
Then Eq. (1.1) has a closed solution on $[t_0,T]$.
}

Proof. By Lemma 2.4 it follows from $12^0)-14^0)$ that $\eta_{N_T,T}(t), \ph t \in [t_0,T]$ is a solution of the inequality (2.3) on $[t_0,T]$ and it follows from the conditions, $12^0), \ph 13^0), \ph 15^0), \ph 16^0)$ that $\zeta_{N_T,T}(t), \ph t \in [t_0,T]$ is a solution of the inequality (2.4) on $[t_0,T]$. It follows from the condition $17^0)$ that $\eta_{N_T,T}(t_0) \ge \eta_{N_T,T}(T)$, and it follows form the condition $18^0)$ that $\zeta_{N_T,T}(t_0) \le \zeta_{N_T,T}(T)$. Therefore, by virtue of  Theorem 2.1  Eq. (1.1) has a closed solution on $[t_0,T]$. The theorem is proved.

\vsk

{\bf Remark 5.1.} {\it The conditions $14^0), \ph 15^0)$ of Theorem 5.5 for $n$ odd are satisfied if, in particular, $d_2(t) = d_{n-1}(t) \equiv 0, \ph d_3(t) > 0, \ph d_{n-2}(t) > 0, \ph d_k(t) = \alpha_k(t) + \beta_k(t), \alpha_k(t) >~ 0, \linebreak \beta_k(t)> 0, \ph k = 5, \ph 7, \ldots, n-4,$
$$
d_4^2(t) - 4 \alpha_5(t) d_3(t) \le 0,
$$
$$
d_6^2(t) - 4\alpha_7(t)\beta_5(t) \le 0,
$$
$$
d_8^2(t) - 4\alpha_9(t)\beta_7(t) \le 0,
$$
$$
\ldots\ldots\ldots\ldots\ldots\ldots\ldots\ldots
$$
$$
d_{n-5}^2(t) - 4\alpha_{n-4}(t)\beta_{n-6}(t) \le 0,
$$
$$
d_{n-3}^2(t) - 4d_{n-2}(t)\beta_{n-4}(t) \le 0, \ph t \in [t_0,T] \ph (n\ge 7)
$$
(since under the above restrictions the "square trinomials"  $\alpha_5(t) u^2 \pm d_4(t)u + d_3(t), \ldots,\linebreak  d_{n-2}(t)u^2 \pm d_{n-3}(t) u +\beta_{n-4}(t)$ are nonnegative for all $t \in [t_0,T], \ph u \in \mathbb{R}$).
Note that the conditions $17^0)$ and $18^0)$ are satisfied if, in particular, $ \il{t_0}{T}\exp\biggl\{\il{t_0}{\tau}a_1(s) d s\biggr\} a_0(\tau) d \tau = 0, \linebreak   \il{t_0}{T}a_1(t) d t \ge 0,$. Indeed, under these restrictions the left part of $17^0)$ is non positive and the left part of $18^0)$ is nonnegative.
}

\vsk

{\bf Example 5.2} {\it For $n=7, \ph a_7(t) = \sin^2 t, \ph  a_6(t) = \sin^2 t \cos t, \ph a_5(t) = 7 \sin^2 t \cos 3 t +~ 2, \linebreak a_4(t) = 4\sin^2 t \arctan t + \sin (\cos t), a_3(t) = 10 \sin^4 t \cos e^t + 2, \ph a_2(t) = \sin^8 t \cos^9 t, \ph t \ge~ t_0, \linebreak \il{t_0}{T}a_1(t) d t \ge 0, \ph \il{t_0}{T}\exp\biggl\{\il{t_0}{\tau}a_1(s) d s\biggr\}a_0(\tau) d\tau = 0$
the conditions of Theorem 5.5 are satisfied.
}

\vsk

{\bf Theorem 5.6.} {\it Let the following conditions be satisfied

\noindent
$19^0)$ for some

$c^+ \ge\max\limits_{t\in[t_0,T]}\il{t_0}{t}\exp\biggl\{-\il{t_0}{\tau}\alpha(s) d s \biggr\}a_0(\tau) d\tau, \ph c^-\ge - \min\limits_{t\in[t_0,T]}\il{t_0}{t}\exp\biggl\{-\il{t_0}{\tau}\alpha(s) d s\biggr\}a_0(\tau) d\tau$\\ the inequalities
$$
\theta_{c^+}(t) \le 1, \phh \theta_{c^-}^-(t)| \le 1, \phh t \in [t_0,T]
$$
are valid,

\noindent
$20^0) \ph c^+\biggl(1 - \exp\biggl\{-\il{t_0}{\tau}\alpha(\tau) d \tau\biggr\}\biggr) \le \il{t_0}{t}\exp\biggl\{-\il{t_0}{\tau}\alpha(s) d s\biggr\}a_0(\tau) d\tau, \ph$

\phantom{aaaaaaaaaaaaaaaaaaaaaaaaa}$c^-\biggl(1 - \exp\biggl\{-\il{t_0}{\tau}\alpha(\tau) d \tau\biggr\}\biggr) \ge \il{t_0}{t}\exp\biggl\{-\il{t_0}{\tau}\alpha(s) d s\biggr\}a_0(\tau) d\tau, \ph$

\noindent
Then Eq. (1.1) has a closed solution $y_*(t)$ on $[t_0,T]$ such that
$$
\theta_{c_-}^-(t) \le y_*(t) \le \theta_{c^+}(t), \phh t \in [t_0,T],
$$
and if $\theta_{c_-}^-(t_0) \le y_*(t_0) \ph (y_*(t_0) \le \theta_{c^+}(t_0))$, then
$$
\theta_{c_-}^-(t) \le y_*(t) \ph (y_*(t) \le \theta_{c^+}(t), \phh t \in [t_0,T].
$$
}

Proof. By Lemma 2.6 it follows from the condition $19^0)$ that $\theta_{z^+}(t)$ and $\theta_{c^-}^-(t)$ are solutions of the inequalities (2.3)   and (2.4) respectively on $[t_0,T]$. It is not difficult to verify that the conditions $20^0)$ imply that
$$
\theta_{c^+}(t_0) \ge \theta_{c^+}(T), \phh \theta_{c^-}^-(t_0) \le \theta_{c^-}^-(T)
$$
By Lemmas 2.1 and 2.2 it follows from here that the solutions $y_1(t)$ and $y_2(t)$ of Eq. (1.1) with $y_1(t_0) = \theta_{c^-}^-(t_0), \ph y_2(t_0) = \theta_{c^+}(t_0)$ exist on $[t_0,T]$ and
$$
y_1(t_0) \le y_1(T), \phh y_2(t_0) \ge y_2(T).
$$
By  Theorem 2.1  it follows from here that Eq. (1.1) has a closed solution $y_*(t)$ on $[t_0,T]$ such that
$$
\theta_{c_-}^-(t) \le y_*(t) \le \theta_{c^+}(t), \phh t \in [t_0,T],
$$
and if $\theta_{c_-}^-(t_0) \le y_*(t_0) \ph (y_*(t_0) \le \theta_{c^+}(t_0))$, then
$$
\theta_{c_-}^-(t) \le y_*(t) \ph (y_*(t) \le \theta_{c^+}(t), \phh t \in [t_0,T].
$$
The theorem is proved.

\vsk

{\bf Example 5.3.} {\it Assume $n\ge 2$ is odd, $b_0(t) =\ldots = b_{n}(t) = \nu(t), \ph e_k(t) = (-1)^kb_k(t), \linebreak k=\overline{0,n}, \ph t \in [t_0,T].$ Then, obviously, $y_1(t)\equiv -1, \ph t \in [t_0,T]$ is a solution of Eq. (2.1) on $[t_0,T]$ and $y_2(t)\equiv 1, \ph t \in [t_0,T]$ is a solution of Eq. (3.9) on $[t_0,T]$. It is not difficult to verify that the conditions $(III)$ and $(IV)$ of Theorem 3.2 are satisfied provided $\sum\limits_{k=0}^n(-1)^k a_k(t) \le 0, \ph \sum\limits_{k=0}^n a_k(t) \ge 0, \ph t \in [t_0,T]$.
 Then according to Theorem 3.2 the solutions $y_1^*(t)$ and $y_2^*(t)$  of Eq. (1.1) with $y_1^*(t_0)=-1, \ph y_2^*(t_0)=1$ exist on $[t_0,T]$ and $y_1^*(T)\ge y_1^*(t_0), \ph y_2^*(T)\le y_1^*(t_0)$. By  Theorem 2.1 it follows from here that Eq. (1.1) has a closed solution on $[t_0,T]$. Assume $n\ge 2$ is even, $b_0(t) = b_2(t) = \ldots = b_n(t) = \nu(t), \ph b_1(t) = 2 \nu(t), \ph e_k(t) = (-1)^kb_k(t), \ph k=\overline{0,n}, \ph t \in [t_0,T].$  Then, as for the previous case, $y_1(t)\equiv -1, \ph t \in [t_0,T]$ is a solution of Eq. (2.1) on $[t_0,T]$ and $y_2(t)\equiv 1, \ph t \in [t_0,T]$ is a solution of Eq. (3.9) on $[t_0,T]$. In this case we can also use Theorem 3.2 to conclude that Eq. (1.1) has a closed solution on $[t_0,T]$ provided $\sum\limits_{k=0}^n(-1)^k a_k(t) \le 0, \ph  \sum\limits_{k=0}^n a_k(t) \ge 0, \ph   t \in [t_0,T].$
Thus, for $n \ge 2$ Eq. (1.1) has a closed solution provided
$$
\sum\limits_{k=0}^n(-1)^k a_k(t) \le 0, \ph  \sum\limits_{k=0}^n a_k(t) \ge 0, \phh   t \in [t_0,T].
$$
}

\vsk

\centerline{\bf References}

\vsk

\noin
1.  G. A. Grigorian,  On two comparison tests for second-order linear ordinary differential \linebreak \phantom{a}
equations, Diff. Urav., vol 47, (2011), 1225--1240 (in Russian), Diff. Eq., vol. 47, (2011),   \linebreak \phantom{a}
1237--1252  (in English).

\noin
2.  G. A. Grigorian,   Two comparison criteria for scalar Riccati equations and some of  \linebreak \phantom{a} their applications
Izv. Vissh. Uchebn. Zaved. Mat. (2012), no 11, 20-25.
Russian \linebreak \phantom{a}  Math. (Iz. VUZ) 56 (2.12), no. 11, 17-30.

\noin
3.   G. A. Grigorian,   Interval oscillation criteria for linear matrix Hamiltonian systems.  \linebreak \phantom{a}
Rocky Mountain J. Math. 50 (2000), no. 6, 2047-2057.

\noin
4.  G. A. Grigorian, New reducibility criteria for systems of two linear first-order ordinary  \linebreak \phantom{a} differential equations.
Monatsh. Math. 198 (2022), no. 2. 311-322.

\noin
5.  G. A. Grigorian, Oscillation and non-oscillation criteria for linear nonhomogeneous  \linebreak \phantom{a} systems of two first-order ordinary differential equations.
J. Math. Anal. Appl. 507  \linebreak \phantom{a} (2022), no. 1, Paper No. 125734, 10 pp.

\noindent
6.  G. A. Grigorian, Oscillatory and non-oscillatory criteria for linear four-dimensional   \linebreak \phantom{a} Hamiltonian systems
Math. Bohem. 146 (2021), no. 3, 289-304.

\noin
7.  G. A. Grigorian,  On the reducibility of systems of two linear first-order ordinary  \linebreak \phantom{a} differential equations.
Monatsh. Math. 195 (2021), no. 1, 107-117.

\noin
8.  G. A. Grigorian, Oscillation criteria for linear matrix Hamiltonian systems.
Proc. \linebreak \phantom{a} Amer. Math. Soc. 148 (2020), no. 8, 3407-3415.

\noin
9.  G. A. Grigorian,  Oscillatory and non oscillatory criteria for the systems of two linear  \linebreak \phantom{a} first order two by two dimensional matrix ordinary differential equations.
Arch. Math.\linebreak \phantom{a} (Brno) 54 (2018), no. 4, 169-203.

\noin
10.  G. A. Grigorian,   Stability criterion for systems of two first-order linear ordinary   \linebreak \phantom{a}  differential equations.
Mat. Zametki 103 (2018), no. 6, 831-840.
Math. Notes 103 \linebreak \phantom{a} (2018), no. 5-6, 892-900.

\noin
11. G. A. Grigorian,  Oscillatory criteria for the systems of two first-order linear differ-\linebreak \phantom{a} ential  equations.
Rocky Mountain J. Math. 47 (2017), no. 5. 1497-1524.

\noin
12.  G. A. Grigorian, Some properties of the solutions of third order linear ordinary   \linebreak \phantom{a} differential equations
Rocky Mountain J. Math. 46 (2016), no. 1, 147-168.

\noin
13.  G. A. Grigorian, Stability criteria for systems of two first-order linear ordinary \linebreak \phantom{a} differential  equations.
Math. Slovaca 72 (2022), no. 4, 935-944.

\noin
14. Yu. Ilyashenko, Centennial History of Hilbert's 16th Problem.
Bul. (New Series) of   \linebreak \phantom{a} the
AMS, Vol. 39, Num. 3, pp. 301--354.

\noin
15. M. Briskin and Y. Yumdin, Tangential version of Hilbert 16th problem for the Abel \linebreak \phantom{a}  equation. Moscow Mathematical Journal, vol. 5, No. 1, 2005, 23--53.

\noin
16. A. Gassull, J. Libre, Limit cycles for a class of Abel equations. SIAM J. Math. Anal. \linebreak \phantom{a} 21 (1990), 1235--1244.

\noin
17. A. Gassull, R. Prohens, J. Torregrossa, Limit cycles for  rigid cubic systems. J. Math.\linebreak \phantom{a} Anal. Appl.  303 (2005), 391--404.

\noin
18. P. J. Torres, Existence of closed solutions for a polynomial first order differential \linebreak \phantom{a}  equation. J. Math. Anal. Appll. 328 (2007) 1108--1116.

\noin
19. B. P. Demidovich, Lectures on the mathematical theory of stability. Moskow, \linebreak \phantom{a} "Nauka", 1967.

\noin
20. S. E. Ariaku, E. C. Mbuh, C. C. Asogva, P. U. Nwokoro. Lower and Upper Solu-\linebreak \phantom{a} tions of First Order Non-Linear Ordinary Differential Equations. International Journal \linebreak \phantom{a} of Scientific Engineering and Science, Vol. 3, Issue 11, 2019, pp. 59--61.

\noin
21. R. E. Edwards, A formal background to mathematics (Springer-Verlag, New York, \linebreak \phantom{a} Haidelberg, Berlin, 1980).

\end{document}